\documentclass[12pt,reqno,a4paper]{article}
\usepackage{geometry}
\usepackage[blocks]{authblk}
\usepackage{amsmath,amsthm,amstext,amscd,amssymb,euscript,url,mathrsfs,euscript,dsfont,calrsfs,mathtools,comment,relsize,xcolor}
\usepackage{hyperref}
\usepackage{epsfig}
\usepackage[inline]{enumitem}
\usepackage[absolute,overlay]{textpos}
\usepackage{graphicx}
\usepackage{subfig}
\usepackage{tikz}
\usetikzlibrary{matrix}
\usetikzlibrary{arrows,positioning}

\renewcommand{\P}{\mathbb{P}}

\newcommand{\R}{\mathbb R}
\newcommand{\N}{{\mathbb N}}

\newcommand{\E}{\mathbb E}

\renewcommand{\Pr}{\mathbb P}
\newcommand\var{\textup{Var}}
\newcommand\cov{\textup{Cov}}
\renewcommand{\phi}{\varphi}
\newcommand{\epsi}{\ensuremath{\epsilon}}

\newcommand{\si}{\ensuremath{\sigma}}

\def\1{{\mathchoice {\rm 1\mskip-4mu l} {\rm 1\mskip-4mu l}
{\rm 1\mskip-4.5mu l} {\rm 1\mskip-5mu l}}}

\newtheorem{theorem}{{\small T}{\scriptsize HEOREM}}[section]
\newtheorem{corollary}{{\bf{\small C}{\scriptsize OROLLARY}}}[section]
\newtheorem{proposition}{{\bf{\small P}{\scriptsize ROPOSITION}}}[section]
\newtheorem{lemma}{{\bf{\small L}{\scriptsize EMMA}}}[section]
\newtheorem{remark}{{\bf{\small R}{\scriptsize EMARK}}}[section]
\newtheorem{definition}{{\bf{\small D}{\scriptsize EFINITION}}}[section]
\newtheorem{induction}{{\bf{\small I}{\scriptsize NDUCTIVE HYPOTHESIS}}}[section]

\renewenvironment{proof}[1]
{\noindent{{\bf{\small{ P}{\scriptsize ROOF}}}.}\hspace{0.1cm} #1} {$\;\qed$\newline}

\newcommand{\beq}{\begin{eqnarray}}
\newcommand{\eeq}{\end{eqnarray}}

\newcommand{\ba}{\begin{align*}}
\newcommand{\ea}{\end{align*}}

\newcommand{\be}{\begin{equation}}
\newcommand{\ee}{\end{equation}}

\newcommand{\bl}{\begin{lemma}}
\newcommand{\el}{\end{lemma}}

\newcommand{\br}{\begin{remark}}
\newcommand{\er}{\end{remark}}

\newcommand{\bt}{\begin{theorem}}
\newcommand{\et}{\end{theorem}}

\newcommand{\bd}{\begin{definition}}
\newcommand{\ed}{\end{definition}}

\newcommand{\bind}{\begin{induction}}
\newcommand{\eind}{\end{induction}}

\newcommand{\bp}{\begin{proposition}}
\newcommand{\ep}{\end{proposition}}

\newcommand{\bc}{\begin{corollary}}
\newcommand{\ec}{\end{corollary}}

\newcommand{\bpr}{\begin{proof}}
\newcommand{\epr}{\end{proof}}

\newcommand{\bi}{\begin{itemize}}
\newcommand{\ei}{\end{itemize}}

\newcommand{\ben}{\begin{enumerate}}
\newcommand{\een}{\end{enumerate}}

\newcommand{\caF}{{\mathcal F}}

\newcommand{\caI}{{\mathcal I}}

\newmuskip\pFqmuskip

\newcommand\pFq[6][8]{%
	\begingroup 
	\pFqmuskip=#1mu\relax
	\mathcode`\,=\string"8000
	\begingroup\lccode`\~=`\,
	\lowercase{\endgroup\let~}\pFqcomma
	{}_{#2}F_{#3}{\left[\genfrac..{0pt}{}{#4}{#5};#6\right]}%
	\endgroup
}
\newcommand{\pFqcomma}{\mskip\pFqmuskip}

\newcommand{\norm}[1]{\left\lVert#1\right\rVert}

\newtheorem{asu}{Assumption}
\newcommand{\basu}{\begin{asu}}
\newcommand{\easu}{\end{asu}}

\begin{document}

\title{Ergodic properties of the harmonic process}
\author[1]{Frank Redig
\thanks{Email: \texttt{F.H.J.Redig@tudelft.nl}}
}
\author[2]{Berend van Tol
\thanks{Email: \texttt{B.T.vanTol@tudelft.nl}}
}

\affil[1,2]{Institute of Applied Mathematics, Delft University of Technology, Delft,
The Netherlands}

\date{\today}

\maketitle
\begin{abstract}
In this paper we study detailed fluctuation results for a class of non-equilibrium steady states. The main example is the boundary driven harmonic model \cite{frassek2022exact}.
In this model, the non-equilibrium  steady state (NESS) is a mixture of products of geometric distributions, of which the local parameters are in turn distributed as uniform order statistics. For such a NESS, we prove law of large numbers,  central limit theorem and large deviation results for fields of general local functions (generalizing the density field). We also obtain quantitative results on the deviation from local equilibrium. 
    
\end{abstract}

\section{Introduction}
Understanding the microscopic and emergent macroscopic properties of non-equilibrium steady states (NESS) is an important challenge in mathematical statistical physics. At present only for a limited class of boundary driven models one can access microscopic details of the NESS such as the expected profile, correlation functions and large deviation rate functions. Important examples include boundary driven stochastic systems such as the (a)symmetric exclusion process \cite{derrida1998exactly} which is exactly solvable via the matrix ansatz solution. More recently, the harmonic model and a more general class of exactly solvable models were introduced in \cite{frassek2022exact,frassek2020non}.
In these boundary-driven generalized harmonic models, the equilibrium measures are product measures, and in the presence of boundary reservoirs the NESS can be written as a mixture of product measures with equilibrium marginals. In the simplest setting of the harmonic process on a chain $\{1,\ldots, N\}$ with left (resp.\ right) boundary reservoir parameter $\theta_L$ (resp. $\theta_R>\theta_L$), this NESS is a product of geometric distributions, where the parameters of the geometric distributions are themselves distributed as order statistics of independent uniforms on the interval $[\theta_L,
\theta_R]$, see \cite{carinci2024solvable,giardina2025intertwining,carinci2025large}. Therefore, in that setting the measure describing the mixture is explicitly given as the joint distribution of those uniform order statistics. The fact that the NESS can be written as a mixture of product measures with equilibrium marginals holds for a much larger class of models on general graphs, including, e.g., the boundary driven KMP model, see \cite{giardina2025intertwining,
de2024hidden}. However, at present, only in the boundary driven generalized harmonic models on a chain with reservoirs at the left and right ends, the measure which describes this mixture is known in closed form. As a consequence, in those models many detailed properties of the corresponding NESS can be derived, including a large deviation principle for the density profile \cite{carinci2025large} with the rate function predicted from macroscopic fluctuation theory \cite{bertini2, bertini2007stochastic}.

In this paper we revisit the boundary driven harmonic model,
and from the representation of the NESS as a mixture of product states, we derive several detailed results about fluctuations of the density field and more general fields associated to local functions. Roughly speaking, these results follow from a combination of precise fluctuation results for product measures with slowly varying parameters, combined with precise fluctuation results for order statistics. With minor modifications, these results hold for the entire one-parameter family of boundary driven
(generalized) harmonic models and exactly solvable heat conduction models \cite{franceschini2023integrable}. Indeed, for all these models, the NESS can be described as a mixture of product measures, where the measure which describes the mixtures is expressed in terms of generalized order statistics (Dirichlet process). 

In summary, in our paper we obtain the following main results for the NESS of the boundary driven harmonic model.
\begin{enumerate}
\item Strong law of large numbers for the profile of a local function.
\item Central limit theorem for the profile of a local function.
\item Local equilibrium and 
deviations from local equilibrium for finite system size.
\item Large deviations for the empirical profile of a local function.
\end{enumerate}
The remainder of our paper is organized as follows.
First, we define the precise mixture of product states which we study (Section 2). Second, we prove fluctuation results for order statistics, i.e., the random parameters (Section 3).
We then use these fluctuation results to prove the results announced above.
More precisely, in Section 4 we prove local equilibrium (in the large system-size limit $N\to\infty$) and quantitative deviations from it for finite system size $N$, in Section 5 we prove the law of large numbers for fields of general local functions, in Section 6 we prove central limit theorems for the fluctuation fields of general local functions, and finally in Section 7 we prove large deviations for
 fields of general (bounded) local functions.

\section{Mixed product states} \label{sec:mixture_structure}
In this section we describe the mixed product states which we study, where we focus on mixtures of geometric distributions, inspired by the stationary measure of the boundary driven harmonic model 
\cite{frassek2022exact}.
We consider particle systems on a finite chain
$\{1,\ldots, N\}$ with left and right boundary reservoirs. The configuration space of such a particle system is given by
$\Omega_N= \N^{\{1,\ldots,N\}}$ where $\N= \{0,1,2,\ldots\}$. For $\eta\in\Omega_N$, and $i\in \{1,\ldots,N\}$, $\eta_{i,N}$ represents the number of particles at lattice site $i$.
Let $\nu_\theta$ denote a probability measure on $\N$, parametrized by a parameter $\theta\in [0,\infty)$. The main example considered in this paper is 
the geometric distribution, which we parametrize as follows:
\be\label{geodi}
    \nu_\theta(n) = \frac{1}{1+\theta} \Big( \frac{\theta}{1+\theta} \Big)^n .
\ee
Notice that in this parametrization, the parameter $\theta$  coincides with the expectation, i.e., 
\[
\sum_n n\nu_\theta(n)= \theta.
\]

Then we can define associated product measures on the configuration space $\Omega_N$ given by
$\otimes_{i=1}^N \nu_{\theta(i)}$.
Here $\theta: \{1,\ldots, N\}\to [0,\infty)$ denotes what we call the profile of local parameters.
We call a measure $\mu$ a mixture of product measures with geometric marginals if it is of the form

\be\label{mixprod1}
\mu(\eta) = \mathbb{E}\left(\otimes_{i=1}^N \nu_{\Theta_{i,N}}[\eta_{i,N}]\right).
\ee

Here $\Theta_{i,N}, i=1,\ldots, N$ is a collection of non-negative random variables which we call the random (local) parameters.
The expectation $\mathbb{E} $ in \eqref{mixprod1} is over these random parameters 
$\Theta_{i,N}, i=1,\ldots, N$.
Equivalently, $\mu$ is a mixture of product measures of the form $\nu_\theta$ if it is of the form
\be\label{mixprod2}
\mu(\eta) = \int \dots \int \otimes_{i=1}^N \nu_{\theta_{i,N}}[\eta_{i,N}]\Lambda (d\theta_{1,N}\ldots d\theta_{N,N}) 
\ee
where $\Lambda (d\theta_{1,N}\ldots d\theta_{N,N})$ is a probability measure on
(a subset of) $[0,\infty)^N$, which corresponds in the notation \eqref{mixprod1} to the joint distribution of the random parameters $\Theta_{i,N}, i=1,\ldots, N$. This measure
$\Lambda (d\theta_{1,N}\ldots d\theta_{N,N})$
is the measure which describes the mixture, and is called the ``mixture measure''. 

Motivated by the NESS in the harmonic model, \cite{carinci2024solvable, carinci2025large}, in this paper we restrict to mixtures of products of geometric distributions, where the parameters are distributed as order statistics. That is, the random parameters  $\{\Theta_{i,N}\}_{1 \leq i \leq N}$ are defined as follows. Let $\theta_L, \theta_R \in \R^+$ be such that $\theta_L<\theta_R$. These parameters, $\theta_L$ and $ \theta_R$, correspond to the boundary parameters of the Harmonic models. We define $N$ independent, uniformly distributed random variables on $[\theta_L, \theta_R]$:

\be
U^{\theta_L,\theta_R}_{1}, U^{\theta_L,\theta_R}_{2} \dots U^{\theta_L,\theta_R}_{N} \sim \text{Unif}([\theta_L,\theta_R]).
\ee

The parameters $\Theta_{1,N},..,\Theta_{N,N}$ are then the order statistics of these $N$ uniform random variables $\{U^{\theta_L, \theta_R}_i\}_{1\leq i \leq N}$, i.e.,  for $i\in \{1,...,N\}$,
\be\label{ordthet}
    \Theta_{i,N} := U^{\theta_L,\theta_R}_{i:N},
\ee

{\color{black}
Whenever endpoint indices $i=0$ or $i=N+1$ occur, we use the convention
$\Theta_{0,N}:=\theta_L$ and $\Theta_{N+1,N}:=\theta_R$.
}

with $U^{\theta_L,\theta_R}_{i:N}$ the $i$-th order statistic of $U^{\theta_L, \theta_R}_1, U^{\theta_L, \theta_R}_2, \dots, U^{\theta_L, \theta_R}_N$. Notice that one can construct $\{U^{\theta_L,\theta_R}_{i:N}\}_{1 \leq i\leq N}$ from independent standard uniform random variables $U_1, U_2, \dots, U_N$ via

\be \label{eq:Theta_in_terms_stdU}
U^{\theta_L,\theta_R}_{i:N} := \theta_L + (\theta_R-\theta_L) U_{i:N}
\ee

with $U_1,...,U_N \sim U([0,1])$ and $U_{i:N}$ the $i$-th order statistic. 

We further denote by $\rho_{\theta_L, \theta_R}: [0,1]\to [0,\infty)$ the function
\be
\rho_{\theta_L, \theta_R}(x) := [\theta_L + (\theta_R-\theta_L) x] \mathds{1}(0 \leq x \leq 1).
\ee

\begin{remark}
In what follows, we will often prove results about the mixed measure \eqref{mixprod2} and the
corresponding random variable $\eta_{i,N}, \Theta_{i,N}, i=1, \ldots, N$.
We think of these random variables as defined on a common underlying probability space
(jointly for all values of $N$). This implies that we can speak about almost sure asymptotic results
(as $N\to\infty$).
\end{remark}

\section{Fluctuation properties of the random parameters}
We start this section with a basic result on moments of uniform order statistics, 
which are known in closed form (see e.g. equation 2.1 in \cite{mathai2002product}). Proposition \ref{prop:moments_order_statistics} below states this explicit expression, which we will use occasionally throughout this text. 
\bp \label{prop:moments_order_statistics}
{\color{black}Let $\alpha_1,\ldots,\alpha_N\in\R$ satisfy
$\alpha_1+\cdots+\alpha_j+j>0$ for every $j\in\{1,\ldots,N\}$.} We have

\be \label{eq:exp_products_order_stat}
\E(U_{1:N}^{\alpha_1} \cdot U_{2:N}^{\alpha_2}  \dots U_{N:N}^{\alpha_N}) = \Gamma(N+1) \left \{\prod_{j=1}^N \frac{\Gamma(\alpha_1+ \dots \alpha_j+j)}{\Gamma(\alpha_1+ \dots \alpha_j+j+1)} \right\}.
\ee

In particular, we have the following expression for the $k$-th moment, $k \in \N$, of $U_{r:N}$, $r \in  \{1,\dots,N\}$.

\be \label{eq:moments_order_stat}
\E(U_{r:N}^k) = \frac{r(r+1) \dots (r+k-1)}{(N+1)(N+2) \dots (N+k)}.
\ee
\ep

In the rest of this section we state and prove several ergodic results for the order statistics $\{\Theta_{i,N}\}_{1 \leq i \leq N}$, which we call from now on ``the random parameters'' (in the spirit of section 2). More precisely, in Proposition \ref{prop:lln_theta} we prove a strong law of large numbers (SLLN) and in Proposition \ref{prop:clt_Theta} the central limit theorem (CLT). These results are essentially the SLLN and CLT in chapter 19 of \cite{shorack2009empirical} extended to hold for local functions of the random parameters. More precisely, the SLLN and CLT in \cite{shorack2009empirical} hold for general L-statistics, i.e. random variables of the form
\be
\sum_{i=1}^N h(\Theta_{i,N}) \phi\left(\frac{i}{N}\right),
\ee
for some appropriate functions $\phi$ and $h$ on $\R$. In this paper we are interested in sums of local functions of the random parameters:
{\color{black}
\begin{align}
\sum_{i=0}^{N-k} (\tau_i h)[\Theta] \phi \left(\frac{i}{N}\right) &= \sum_{i=0}^{N-k} h(\Theta_{i+1,N}, \dots, \Theta_{i+k,N}) \phi \left(\frac{i}{N}\right)\\
&= h(\Theta_{1,N}, \dots, \Theta_{k,N}) \phi(0) + \dots + h(\Theta_{N-k+1,N}, \dots, \Theta_{N,N}) \phi \left(\frac{N-k}{N}\right). \nonumber
\end{align}
}

In later sections we use these results to prove the corresponding statements for the random variables 
$\{\eta_{i,N}\}_{1\leq i\leq N}$ with distribution given in \eqref{mixprod1}.\\

Here and in the rest of this paper we define the operator $\tau_i$ acting on local functions $f(x) = f(x_1, \dots, x_k)$ as
$(\tau_i f)[x]:= f(x_{i+1}, \dots, x_{i+k})$. Let $E$ be either $E = \N$ or $E = \R$. Throughout this text we call for fixed $k\in \N$ the function $f:E^N\rightarrow \R$, with $N\geq k$, a local function with dependence set $\{1,2, \ldots,k\}$ if $f(x)$ only depends on $x_1, \ldots, x_k$. That is, with slight abuse of notation, $f(x) = f(x_1, x_2, \ldots, x_k)$. %

In subsection \ref{subsec:concentration_Theta} we state and prove a concentration inequality for the random parameters which will later be used to prove local equilibrium for $\{\eta_{i,N}\}_{1 \leq i \leq N}$. Finally we state a large deviation theorem for the random parameters in subsection \ref{subsec:LDP_Theta}, from which we will later again find corresponding results for the $\{\eta_{i,N}\}_{1 \leq i \leq N}$ random variables.

\subsection{Law of large numbers for the random parameters} 

\bp[SLLN for the random parameters]\label{prop:lln_theta}
Let $\phi: [0,1] \rightarrow \R$ be differentiable with bounded {\color{black}first and second} derivative and let $h: [\theta_L,\theta_R]^N \rightarrow \R$ be both continuous and local with dependence set $\{1,\dots,k\} \subset \N$, $k\in \N$. Then,

\be
\frac{1}{N} \sum_{i=0}^{N-k} (\tau_i h)[\Theta] \phi \Big(\frac{i}{N} \Big) \rightarrow \int_0^1 h(\rho_{\theta_L,\theta_R}(x) ) \phi(x) dx
\ee

almost surely as $N \rightarrow \infty$. Here it is understood that for a real number $\rho \in [\theta_L, \theta_R]$ we write $h(\rho)$ to denote 

\be
{\color{black}h(\rho, \rho, \ldots, \rho).}
\ee

That is, $h(\rho)$ is the value of $h$ acting on $(\rho, \rho, \ldots, \rho) \in {\color{black}[\theta_L, \theta_R]^N}$.
 
\ep

\bpr
It is sufficient to show the result for the case where $h$ is a monomial of the form $h(x_1, \dots, x_k) = x_1^{p_1} \cdot x_2^{p_2} \dots x_k^{p_k}$, $p_1, \dots, p_k \in \N$, since the polynomials are dense in $C([\theta_L,\theta_R]^{k})$. Assume wlog that for all $x \in [0,1]:$ $\phi(x) \geq 0$.\\

Theorem 3 in Chapter 19 of \cite{shorack2009empirical} gives that for $p = p_1 + p_2 + \dots + p_k \in \N$

\be \label{eq:shorack_SLLN}
\frac{1}{N} \sum_{i=1}^N (\Theta_{i,N})^p \phi \left(\frac{i}{N} \right) \rightarrow \int_0^1 (\rho_{\theta_L, \theta_R}(x))^p \phi(x) dx
\ee

Furthermore, by the definition of order statistics, we have that the random parameters $\{\Theta_{i,N}\}_{1\leq i \leq N}$, are monotone in $i$. Hence we almost surely have

\be
\Theta_{1,N}^p \leq h(\Theta) = \Theta_{1,N}^{p_1} \Theta_{2,N}^{p_2} \dots \Theta_{k,N}^{p_k} \leq \Theta_{k,N}^p.
\ee

As a consequence

\be \label{eq:lbd_SLLN_Theta}
\frac{1}{N} \sum_{i=0}^{N-k} (\Theta_{i,N})^{p} \phi \left( \frac{i}{N} \right) \leq \frac{1}{N} \sum_{i=0}^{N-k} (\Theta_{i+1,N}^{p_1} \Theta_{i+2,N}^{p_2} \dots \Theta_{i+k,N}^{p_k}) \phi \left( \frac{i}{N} \right)
\ee

and

\be \label{eq:ubd_SLLN_Theta}
\frac{1}{N} \sum_{i=0}^{N-k} (\Theta_{i+k,N})^{p} \phi \left( \frac{i}{N} \right) \geq \frac{1}{N} \sum_{i=0}^{N-k} (\Theta_{i+1,N}^{p_1} \Theta_{i+2,N}^{p_2} \dots \Theta_{i+k,N}^{p_k}) \phi \left( \frac{i}{N} \right).
\ee

Using \eqref{eq:shorack_SLLN} it is straightforward to see that the left hand side of \eqref{eq:lbd_SLLN_Theta} and the left hand side of \eqref{eq:ubd_SLLN_Theta} both converge almost surely to the same limit:

{\color{black}
\begin{align}
\frac{1}{N} \sum_{i=0}^{N-k} (\Theta_{i+k,N})^{p} \phi \left( \frac{i}{N} \right) &=  \frac{1}{N} \sum_{i=k}^{N} (\Theta_{i,N})^{p} \phi \left( \frac{i-k}{N} \right) \\
&= \frac{1}{N} \sum_{i=k}^{N} (\Theta_{i,N})^{p} \left[\phi \left( \frac{i}{N}\right) - \frac{k}{N}\phi'(\xi_{i,N})\right] \nonumber \\
&= \frac{1}{N} \sum_{i=0}^{N} (\Theta_{i,N})^{p}\phi \left( \frac{i}{N}\right) + O(N^{-1})\xrightarrow{a.s.} \int_0^1 (\rho_{\theta_L, \theta_R}(x))^p \phi(x) dx, \nonumber
\end{align}
where $\xi_{i,N}$ lies between $(i-k)/N$ and $i/N$. The error is
$O(N^{-1})$ because $\phi'$ is bounded.
}

and

\begin{align}
\frac{1}{N} \sum_{i=0}^{N-k} (\Theta_{i,N})^{p} \phi \left( \frac{i}{N} \right) &= \frac{1}{N} \sum_{i=0}^{N} (\Theta_{i,N})^{p}\phi \left( \frac{i}{N}\right) + O(N^{-1})\\
&\xrightarrow{a.s.} \int_0^1 (\rho_{\theta_L, \theta_R}(x))^p \phi(x) dx.\nonumber
\end{align}

We conclude that 
\begin{align}
\frac{1}{N} \sum_{i=0}^{N-k} (\Theta_{i+1,N}^{p_1} \Theta_{i+2,N}^{p_2} \dots \Theta_{i+k,N}^{p_k}) &\phi \left( \frac{i}{N} \right) = \frac{1}{N} \sum_{i=0}^{N-k} (\tau_i h)[\Theta] \phi \left( \frac{i}{N} \right)\\
&\xrightarrow{a.s.} \int_0^1 (\rho_{\theta_L, \theta_R}(x))^p \phi(x) dx = \int_0^1 h(\rho_{\theta_L, \theta_R}(x) \mathds{1}) \phi(x) dx.\nonumber
\end{align}

\epr
\subsection{Central limit theorem for the random parameters}
\bp [CLT for the random parameters] \label{prop:clt_Theta}
Let $\phi:[0,1] \rightarrow \R$ be a continuous function and let $h:[\theta_L,\theta_R]^N \rightarrow \R$ be {\color{black}continuously differentiable} and local with dependence set $\{1, \dots,k\} \subset \N$, $k\in \N$. Then

\be \label{eq:fluki}
 \lim_{N\rightarrow \infty} \frac{1}{\sqrt{N}} \sum_{i=0}^{N-k} \Big((\tau_i h)[\Theta] - \E((\tau_i h)[\Theta])\Big) \phi \Big(\frac{i}{N}\Big) \stackrel{\text{d}}{=} \mathcal{N}(0, \sigma_T^2)
\ee

Here, $\sigma_T^2$ is given by

\be
\sigma_T^2 = (\theta_R-\theta_L)^2 \int_{0}^{1} \int_0^1 [s \wedge t -st] \phi(s) \phi(t) h'(\rho_{\theta_L, \theta_R}(s)) h'(\rho_{\theta_L, \theta_R}(t)) ds dt.
\ee

with $h': [\theta_L,\theta_R] \rightarrow \R$ defined as

\be
h'(\theta) = \left[\frac{\partial} {\partial t} h(t, \dots,t) \right]_{t = \theta}.
\ee
\ep

{\color{black}
\bpr
Define 
$\overline h(\theta):=h(\theta,\ldots,\theta)$.  Since $h$ is continuously differentiable on the compact set $[\theta_L,\theta_R]^k$, there exists a constant $C<\infty$ such that, for every $i\in\{0,\ldots,N-k\}$,
\be \label{eq:clt_diagonal_bound}
\left|
(\tau_i h)[\Theta]-\overline h(\Theta_{i+1,N})
\right|
\leq
C\sum_{j=2}^k
\left(\Theta_{i+j,N}-\Theta_{i+1,N}\right).
\ee
Here we used that the order statistics are increasing.  Moreover, for every $j\in\{2,\ldots,k\}$,
\be \label{eq:clt_telescoping_bound}
\begin{aligned}
\sum_{i=0}^{N-k}
\left(\Theta_{i+j,N}-\Theta_{i+1,N}\right)
&=
\sum_{\ell=N-k+2}^{N-k+j}\Theta_{\ell,N}
-\sum_{\ell=1}^{j-1}\Theta_{\ell,N}\\
&\leq (j-1)(\theta_R-\theta_L).
\end{aligned}
\ee
It follows from \eqref{eq:clt_diagonal_bound} and \eqref{eq:clt_telescoping_bound} that
\be \label{eq:clt_diagonal_error}
\begin{aligned}
&\Bigg|
\frac{1}{\sqrt N}\sum_{i=0}^{N-k}
\Big[
(\tau_i h)[\Theta]-\E((\tau_i h)[\Theta])
-\overline h(\Theta_{i+1,N})
+\E\big(\overline h(\Theta_{i+1,N})\big)
\Big]\phi\left(\frac{i}{N}\right)
\Bigg|\\
&\qquad\leq
\frac{2C\norm{\phi}_\infty}{\sqrt N}
\sum_{j=2}^k(j-1)(\theta_R-\theta_L)
\longrightarrow 0.
\end{aligned}
\ee
Consequently, the terms in the sum on the lhs of \eqref{eq:fluki} can be replaced explicitly by their diagonal counterparts:
\be \label{eq:clt_explicit_replacement}
\begin{aligned}
&\frac{1}{\sqrt N}\sum_{i=0}^{N-k}
\Big((\tau_i h)[\Theta]-\E((\tau_i h)[\Theta])\Big)
\phi\left(\frac{i}{N}\right)\\
&\qquad=
\frac{1}{\sqrt N}\sum_{i=0}^{N-k}
\Big(\overline h(\Theta_{i+1,N})
-\E\big(\overline h(\Theta_{i+1,N})\big)\Big)
\phi\left(\frac{i}{N}\right)+o(1)\\
&\qquad=
\frac{1}{\sqrt N}\sum_{j=1}^{N-k+1}
\Big(\overline h(\Theta_{j,N})
-\E\big(\overline h(\Theta_{j,N})\big)\Big)
\phi\left(\frac{j-1}{N}\right)+o(1).
\end{aligned}
\ee
The error in \eqref{eq:clt_explicit_replacement} converges to zero almost surely.  Since
$\Theta_{j,N}=\theta_L+(\theta_R-\theta_L)U_{j:N}$, the last sum in \eqref{eq:clt_explicit_replacement} is 
\be \label{eq:clt_reduced_l_statistic}
\frac{1}{\sqrt N}\sum_{j=1}^{N-k+1}
\left[
\overline h\!\left(\theta_L+(\theta_R-\theta_L)U_{j:N}\right)
-\E\!\left(
\overline h\!\left(\theta_L+(\theta_R-\theta_L)U_{j:N}\right)
\right)
\right]
\phi\left(\frac{j-1}{N}\right).
\ee
Applying Theorem 1 in Chapter 19 of \cite{shorack2009empirical} to \eqref{eq:clt_reduced_l_statistic}, gives convergence to a centered normal distribution with variance
\be \label{eq:clt_reduced_variance}
(\theta_R-\theta_L)^2
\int_0^1\int_0^1
(s\wedge t-st)\phi(s)\phi(t)
\overline h'(\rho_{\theta_L,\theta_R}(s))
\overline h'(\rho_{\theta_L,\theta_R}(t))
\,ds\,dt.
\ee
Finally, we write $\overline h'(\theta)=h'(\theta)$ as in the definition in the statement.  Thus \eqref{eq:clt_reduced_variance} equals $\sigma_T^2$, and the claim follows from \eqref{eq:clt_explicit_replacement}.
\epr
}

\subsection{Concentration for the random parameters} \label{subsec:concentration_Theta}

In Section \ref{sec:local_equilibrium} we will show that the measures $\mu_N$ satisfy \emph{local equilibrium}. This entails that  $\mu_N$ look locally like a product measure if one takes the limit $N \rightarrow \infty$. This property is mathematically formulated in
Definition \ref{def:le}. The reason local equilibrium holds for $(\mu_N)_{N \geq 1}$ is that the random variables $\Theta_{i,N}$ concentrate sufficiently fast on their expectation as $N \rightarrow \infty$:  

\be
\Theta_{i,N} \approx \E(\Theta_{i,N}) = \theta_L + (\theta_R-\theta_L) \frac{i}{N+1} =\rho_{\theta_L, \theta_R}\left( \frac{i}{N+1}\right).
\ee

Moreover, the random parameters concentrate \emph{uniformly}. This is the content of Proposition \ref{prop:maximal_deviation}. 

\bp
\label{prop:maximal_deviation}
    Let $\epsilon>0$ and define $\Delta_{i,N} := |\Theta_{i,N} - \mathbb{E}[\Theta_{i,N}]|$. 

{\color{black}
Put $d:=|\theta_R-\theta_L|$.  Then
\be
 \mathbb{P}\left(\sup_{i\in\{1,\ldots,N\}}\Delta_{i,N}
 \geq\epsilon\right)
 \leq\frac{3d^4}{(N+1)\epsilon^4}. \label{eq:concentration_bd}
\ee
}
\ep

{\color{black}
\begin{proof}
Write
\[
 \Theta_{i,N}=\theta_L+(\theta_R-\theta_L)U_{i:N},
\]
where $U_{i:N}\sim\operatorname{Beta}(i,N+1-i)$ and
$\E[U_{i:N}]=i/(N+1)$.  We first recall that if
$Y\sim\operatorname{Beta}(\alpha,\beta)$ and
$q:=\alpha+\beta$, then
\begin{align}
 \E\left[(Y-\E[Y])^4\right]
 &=
 \frac{3\alpha\beta\left(
 2(\alpha-\beta)^2+\alpha\beta(q+2)
 \right)}
 {q^4(q+1)(q+2)(q+3)}
 \leq\frac{3}{q^2}. \label{eq:beta_fourth_moment_bound}
\end{align}
For $q\leq1$, the inequality follows from $0\leq Y\leq1$.  For
$q\geq1$, it follows by using
$\alpha\beta\leq q^2/4$ and
$(\alpha-\beta)^2\leq q^2$ in the displayed formula.

For $U_{i:N}$ we have $q=N+1$, independently of $i$.  Therefore
\[
 \E[\Delta_{i,N}^4]
 =d^4\E\left[
 \left|U_{i:N}-\frac{i}{N+1}\right|^4
 \right]
 \leq\frac{3d^4}{(N+1)^2}.
\]
Markov's inequality gives
\begin{align}
 \P\left(\sup_{i\in\{1,\ldots,N\}}\Delta_{i,N}
 \geq\epsilon\right)
 &\leq
 \sum_{i=1}^N\P(\Delta_{i,N}\geq\epsilon)\nonumber\\
 &\leq
 \frac{1}{\epsilon^4}\sum_{i=1}^N\E[\Delta_{i,N}^4]\nonumber\\
 &\leq
 \frac{3Nd^4}{(N+1)^2\epsilon^4}
 \leq
 \frac{3d^4}{(N+1)\epsilon^4},
\end{align}
which proves the claim.
\end{proof}
}

\begin{remark}

{\color{black}
One can choose $\epsilon$ as a function of $N$ such that both $\epsilon(N)$ and the right hand side of \eqref{eq:concentration_bd} converge to zero. Choose for example
$\epsilon(N)=(N+1)^{-1/8}$. Then $\epsilon(N)\rightarrow0$ and
\[
 \frac{3d^4}{(N+1)\epsilon^4(N)}
 =\frac{3d^4}{\sqrt{N+1}}
 \xrightarrow{N\to\infty}0.
\]
}
We will use this fact in the proof of Lemma \ref{lemma:le_theta}.
\end{remark}
    
\subsection{Large deviation principle for the random parameters} \label{subsec:LDP_Theta}
We define the sample paths associated to the random variables $\Theta_{1,N}, \Theta_{2,N}, \dots, \Theta_{N,N}$ as 

\be
\Theta_N(t) = \Theta_{\lfloor (N+1) t \rfloor, N}.
\ee
{\color{black}
At $t=0$ and $t=1$, this definition uses the convention
$\Theta_{0,N}=\theta_L$ and $\Theta_{N+1,N}=\theta_R$.}

In \cite{duffy2011sample} is proven that order statistics, like the random parameters $\Theta_{1,N}, \dots, \Theta_{N,N}$, satisfy a large deviation principle (LDP) and a variational expression for the corresponding rate function is provided. The LDP with rate function for our context is stated in Proposition \ref{prop:LDP_order_stat}. For the proof we refer to \cite{duffy2011sample}.


\bp \label{prop:LDP_order_stat}
Let $D[0, 1]$ denote the space of càdlàg
functions on the unit interval, equipped with Skorohod topology. Let $A_{0,1}(\theta_L,\theta_R) \subset D[0, 1]$ denote the closed set of non-decreasing functions $f : [0, 1] \rightarrow \R$ such that $f(x) \geq \theta_L$ and $f(1) = \theta_R$. Then the sample paths $\Theta_{N}(\cdot)$ satisfy a large deviation principle with rate function

\be\label{bonki}
J(u) = 
\begin{cases}
    -\int_0^1 \log \left( \frac{u'(t)}{\theta_R- \theta_L} \right) dt \qquad &\text{if } u \in A_{0,1}(\theta_L,\theta_R) \text{ is strictly increasing}\\
    \infty \qquad & \text{otherwise}
\end{cases}.
\ee
\ep

\color{black}
\br
Notice that the integral in the rhs of \eqref{bonki}
is well-defined because $u$ is increasing and Skorokhod, therefore the derivative exists almost surely.
\er 
\color{black}
\section{Local equilibrium}
\label{sec:local_equilibrium}

In this section we show that the mixture measures $\mu_N$ satisfy local equilibrium. Before we do so, we first formally introduce the notion of local equilibrium in Definition \ref{def:le}. Intuitively, local equilibrium should be understood as follows. Typically, the interacting particle systems which we consider have homogeneous equilibrium (i.e., reversible) product measures of the form $\otimes_{i=1}^N \lambda_\rho$ with parameter $\rho\in \R_{\geq 0}$. Then the system satisfies local equilibrium if $\mu_N$ looks like a product measure with parameter $\rho(x)$ around each macroscopic point $x$ under the scaling $N\rightarrow \infty$. 

\bd[Local Equilibrium] \label{def:le}
Let $\rho: [0,1] \rightarrow \R_{\geq 0}$ be a non-negative continuous function and let $\{\lambda\}_{\theta\geq 0}$ be a family of probability measures on $\N$ indexed by parameter $\theta \in \R_{\geq0}$. We say that the sequence of probability measures $\mu_N$ on $\N^{N}$, $N\in\N$, is local equilibrium with density profile $\rho$ and marginal $\lambda_{\rho}$ if 

\be \label{eq:le}
\lim_{N\to\infty}\int (\tau_{\lfloor x N \rfloor } g(\eta) )
\ d\mu_N(\eta)={\color{black}\int g(\eta) \otimes_{i=1}^k d\lambda_{\rho(x)}(\eta_i)}
\ee

for every fixed {\color{black}$x\in[0,1)$} and every integrable local function $g: \N^N \rightarrow \R $ with dependence set $\{1,\dots,k\}$.
\ed

In Subsection \ref{subsection:le_concentration_random_parameters} below, we first show how local equilibrium is related to the concentration of the random parameters as proved in Proposition \ref{prop:maximal_deviation}. Next we prove a quantitative version of local equilibrium in terms of the self-duality polynomials associated to Harmonic processes. This is done in Subsection \ref{subsec:le_quantitative}.

\subsection{Local equilibrium from the concentration of the random parameters}
\label{subsection:le_concentration_random_parameters}

If the function $g$ in Definition \ref{def:le} is a monomial, then \eqref{eq:le} reduces to \eqref{eq:le_monomials} in Theorem \ref{thm:le_monomials}. Once one has proven local equilibrium for monomial $g$, it is straightforward to extend the result to general local functions. By conditioning on the random parameters,  the left hand side of \eqref{eq:le_monomials} becomes a sum of terms which are of the form

\be
\E\Bigg[ \prod_{j = 1}^k \frac{1}{1+\Theta_{\lfloor xN \rfloor+j,N}} \left(\frac{\Theta_{\lfloor xN \rfloor+j,N}}{1+\Theta_{\lfloor xN \rfloor+j,N}}\right)^{l_j} \Bigg] 
\ee

The uniform concentration in Proposition \ref{prop:maximal_deviation} allows us to replace the random parameters by their expectation as is stated in Lemma \ref{lemma:le_theta}. The local equilibrium property then follows automatically, since $\mu_N$ becomes a product measure if we replace the random parameters by something deterministic.

\bl \label{lemma:le_theta}
For every fixed {\color{black}$x\in[0,1)$} and every $(l_1,\ldots,l_k)\in\mathbb{N}^k$

    \be
    \label{eq:le_theta}
        \E\Bigg[ \prod_{j = 1}^k \frac{1}{1+\Theta_{\lfloor xN \rfloor+j,N}} \left(\frac{\Theta_{\lfloor xN \rfloor+j,N}}{1+\Theta_{\lfloor xN \rfloor+j,N}}\right)^{l_j} \Bigg] \rightarrow   \left(\frac{1}{1+\rho_{\theta_L,\theta_R}(x)}\right)^k \left(\frac{\rho_{\theta_L,\theta_R}(x)}{1+\rho_{\theta_L,\theta_R}(x)}\right)^{\sum_{j = 1}^k l_j}
    \ee

    as $N \rightarrow \infty$.
\el

\bpr
Let $\Delta_{i,N}$ be as in the statement of Proposition \ref{prop:maximal_deviation}, i.e. $\Delta_{i,N} := |\Theta_{i,N} - \E[\Theta_{i,N}]| = |\Theta_{i,N} - \rho_{\theta_L, \theta_R}(i/(N+1))|$.

{\color{black}
Set $\epsilon(N)=(N+1)^{-1/8}$.  Then
$\epsilon(N)\rightarrow0$, and Proposition
\ref{prop:maximal_deviation} gives
\[
 \P\left(
 \sup_{i\in\{1,\ldots,N\}}\Delta_{i,N}\geq\epsilon(N)
 \right)
 \leq\frac{3d^4}{\sqrt{N+1}}
 \xrightarrow{N\to\infty}0,
\]
where $d=|\theta_R-\theta_L|$.
}
We write the left hand side of (\ref{eq:le_theta}) as

    \begin{align} \label{eq:split_lemma_le}
        &\mathbb{P}\left(\sup_{i \in \{1,...,N\}} \Delta_{i,N} \geq \epsilon(N) \right)\\
        &\qquad \qquad \qquad \cdot \E \left[ \prod_{j = 1}^k \frac{1}{1+\Theta_{j+\lfloor xN \rfloor,N}} \left(\frac{\Theta_{j+\lfloor xN \rfloor,N}}{1+\Theta_{j+\lfloor xN \rfloor,N}}\right)^{l_j}\Big| \sup_{i \in \{1,...,N\}} \Delta_{i,N} \geq \epsilon(N) \right] \nonumber\\
        +& \nonumber \\
        &\mathbb{P}\left(\sup_{i \in \{1,...,N\}} \Delta_{i,N} < \epsilon(N) \right) \nonumber\\
        &\qquad \qquad \qquad \cdot \E \left[ \prod_{j = 1}^k \frac{1}{1+\Theta_{j+\lfloor xN \rfloor,N}} \left(\frac{\Theta_{j+\lfloor xN \rfloor,N}}{1+\Theta_{j+\lfloor xN \rfloor,N}}\right)^{l_j}\Big| \sup_{i \in \{1,...,N\}} \Delta_{i,N} < \epsilon(N) \right]\nonumber.
    \end{align}

    Notice that the expectations in \eqref{eq:split_lemma_le} are bounded because the random parameters take values in $[ \theta_L, \theta_R ]$. By Proposition \ref{prop:maximal_deviation}, the first term vanishes. We find that the second term converges to the right hand side of \eqref{eq:le_theta},

    \begin{align}
        &\E \left[ \prod_{j = 1}^k \frac{1}{1+\Theta_{j+\lfloor xN \rfloor,N}} \left(\frac{\Theta_{j+\lfloor xN \rfloor,N}}{1+\Theta_{j+\lfloor xN \rfloor,N}}\right)^{l_j}\Big| \sup_{i \in \{1,...,N\}} \Delta_{i,N} < \epsilon(N) \right]=\\
        &\qquad {\color{black}\left[ \prod_{j = 1}^k \frac{1}{1+\rho_{\theta_L,\theta_R}\left(\frac{j+\lfloor xN \rfloor}{N+1}\right)} \Bigg( \frac{\rho_{\theta_L,\theta_R}\left(\frac{j+\lfloor xN \rfloor}{N+1}\right)}{1+\rho_{\theta_L,\theta_R}\left(\frac{j+\lfloor xN \rfloor}{N+1}\right)} \Bigg)^{l_j} + O(\epsilon(N))\right]}\nonumber\\&\qquad\rightarrow\prod_{j = 1}^k \frac{1}{1+\rho_{\theta_L,\theta_R}(x)} \Big(\frac{\rho_{\theta_L,\theta_R}(x)}{1+\rho_{\theta_L,\theta_R}(x)}\Big)^{l_j}. \nonumber
    \end{align}
    This concludes the proof 
\epr

\bt \label{thm:le_monomials}
Let $k \in \N$, $p_1, p_2, \dots p_k \in \N$ and {\color{black}$x \in [0,1)$}. Then

\be \label{eq:le_monomials}
\E \left[ \prod_{j=1}^k \eta^{p_j}_{\lfloor xN \rfloor + j, N} \right] \rightarrow \prod_{j=1}^k \E[\eta_{1,1}^{p_j}| \Theta_{1,1} = \rho_{\theta_L, \theta_R} (x)] = \prod_{j = 1}^k \int \xi^{p_j} d\nu_{\rho_{\theta_L, \theta_R}(x)}(\xi)
\ee
 as $N\rightarrow \infty$.
\et

\bpr
First notice that 

\begin{align}
\E \left[ \prod_{j=1}^k \eta^{p_j}_{\lfloor xN \rfloor + j, N} \right] &= \E \left[ \prod_{j=1}^k \E \left[\eta^{p_j}_{\lfloor xN \rfloor + j, N} \Big| \Theta \right]\right]\\
&=  \E \left[ \prod_{j=1}^k \int \xi^{p_j} d\nu_{\Theta_{\lfloor xN \rfloor + j, N}}(\xi)  \right]\nonumber
\end{align}

where we have the following expression for the integral

\be
\int \xi^{p} d\nu_{\theta}(\xi) = \sum_{l=0}^\infty l^{p} \left(\frac{1}{1+\theta}\right) \left(\frac{\theta}{1+\theta}\right)^l.
\ee

As a consequence we have

\begin{align}
\E \left[ \prod_{j=1}^k \eta^{p_j}_{\lfloor xN \rfloor + j, N} \right] &= \sum_{l_1, \dots, l_k = 0 }^\infty \E\left[\prod_{j=1}^k l_j^{p_j} \left(\frac{1}{1+\Theta_{\lfloor xN \rfloor + j, N}}\right) \left(\frac{\Theta_{\lfloor xN \rfloor + j, N}}{1+\Theta_{\lfloor xN \rfloor + j, N}}\right)^{l_j}\right]. \nonumber
\end{align}

By Lemma \ref{lemma:le_theta} we have that the right hand side converges to

\begin{align}
&\sum_{l_1, \dots, l_k = 0 }^\infty {\color{black} \left(\frac{1
}{1+\rho_{\theta_L,\theta_R}(x)}\right)^k } \prod_{j=1}^k l_j^{p_j} \left(\frac{\rho_{\theta_L,\theta_R}(x)}{1+\rho_{\theta_L,\theta_R}(x)}\right)^{\sum_{j=1}^k l_j}\\
&= {\color{black} \left(\frac{1
}{1+\rho_{\theta_L,\theta_R}(x)}\right)^k }\prod_{j=1}^k \sum_{l=0}^\infty l^{p_j} \left(\frac{\rho_{\theta_L,\theta_R}(x)}{1+\rho_{\theta_L,\theta_R}(x)}\right)^{l}= \prod_{j = 1}^k \int \xi^{p_j} d\nu_{\rho_{\theta_L, \theta_R}(x)}(\xi).\nonumber
\end{align}
\epr

\subsection{Quantitative local equilibrium and self-duality polynomials}
\label{subsec:le_quantitative}

If the function $g$ is the \emph{self-duality polynomial} associated to the harmonic process, then we show the following quantitative version of local equilibrium:
\[
\int (\tau_{\lfloor x N \rfloor } g(\eta) )
\ d\mu_N(\eta)- {\color{black}\int g(\eta) \otimes_{i=1}^k d\lambda_{\rho(x)}(\eta_i)}  = O(1/N)
\]
where $\mu_N$ is the mixed product measure in \eqref{mixprod1}. In other words, we show that as far as expectation of self-duality polynomials are concerned, the deviation from local equilibrium is of order $1/N$. This is the content of Theorem \ref{thm:le_quant}. Before we state and prove this theorem, we first introduce the self-duality polynomials.

\bd
The self-duality polynomials $D_N: \N^{\{1,\dots,N\}} \times  \N^{\{1,\dots,N\}} \rightarrow \R$ associated with the Harmonic model are defined by

\be
D_N(\eta,\xi) = \prod_{i=1}^N {\eta_{i,N} \choose \xi_{i,N}} \mathds{1}(\eta_{i,N} \geq \xi_{i,N}).
\ee
\ed

The function $D_N$ has the following property, which is essential in the proof of Theorem \ref{thm:le_quant}. 
\begin{lemma}\label{simplem}
Conditioned on the random parameters the expectation of $D_N(\eta,\xi)$, for a given $\xi$, is given by
\be \label{eq:prop_duality_pol}
\E\big[D_N (\eta, \xi)  | \Theta \big] = \prod_{i=1}^N \Theta_{i,N}^{\xi_{i,N}}.
\ee
\end{lemma}
\begin{proof}
The statement of the lemma follows from straightforward computation. Let $p = \Theta_{1,N}/(\Theta_{1,N}+1)$, $q= 1-p$ and $n=\eta_{1,N}$, then

\begin{align} \label{eq:one_site_duality_expectation}
    &\E\left[{\eta_{1,N} \choose k} \Big| \Theta_{1,N}\right] =\\
    &\sum_{n \geq k} {n \choose k } p^n q = \frac{q p^k}{k!} \sum_{n \geq k} (n-k+1) \dots n \cdot p^{n-k} = \frac{q p^k}{k!} \frac{\partial^k}{\partial p^k} \left(\frac{1}{1-p}\right) = \left(\frac{p}{1-p}\right)^k = \Theta_{1,N}^k.\nonumber
\end{align}

Equation \eqref{eq:prop_duality_pol} now follows, because $\mu_N$ conditioned on $\Theta$ is the product measure $\prod_{i=1}^N \nu_{\Theta_{i,N}}$.
\end{proof}

\bt \label{thm:le_quant}
Let $k \in \N \setminus\{0\}$ and $p_1, \dots, p_k \in \N$. For all $x \in (0,1)$ we have
\be
\E\left[ D_N (\eta, p_1 \delta_{\lfloor xN \rfloor + 1} + \dots + p_k \delta_{\lfloor xN \rfloor + k}) \right]  =  \rho_{\theta_L,\theta_R}(x)^{p_1 + \dots + p_k} +O(N^{-1}).
\ee
\et

\bpr
We first condition on the random parameters and then use the property of the self-duality polynomials,

\begin{align}
\E \Bigg[ D_N \Bigg(\eta, \sum_{j=1}^k p_j \delta_{\lfloor xN \rfloor + j}\Bigg) \Bigg] &= \E \Bigg[ \E\Bigg(D_N \Bigg(\eta, \sum_{j=1}^k p_j \delta_{\lfloor xN \rfloor + j}\Bigg)  \Bigg| \Theta \Bigg)\Bigg] = \E \Bigg[ \prod_{j=1}^k \Theta_{\lfloor xN \rfloor+j,N}^{p_j}  \Bigg].
\end{align}

We now write the random parameters in terms of the standard uniform order statistics as in \eqref{eq:Theta_in_terms_stdU} 

\begin{align} \label{eq:loc_equi_expansion}
&\E \Bigg[ D_N \Bigg(\eta, \sum_{j=1}^k p_j \delta_{\lfloor xN \rfloor + j}\Bigg) \Bigg]\\
&= \E \Bigg[ \prod_{j=1}^k \Theta_{\lfloor xN \rfloor+j,N}^{p_j}  \Bigg]\nonumber\\
&= \E \Bigg[ \prod_{j=1}^k \left(\theta_L + (\theta_R - \theta_L) U_{\lfloor xN \rfloor+j:N} \right)^{p_j}  \Bigg] \nonumber\\
& =\E \Bigg[ \prod_{j=1}^k {\color{black}\sum_{l = 0}^{p_j}} {p_j \choose l} \theta_L^{p_j - l} (\theta_R - \theta_L)^{l} U_{\lfloor xN \rfloor+j:N}^l  \Bigg] \nonumber\\
& = \sum_{l_1 = 0}^{p_1} \dots \sum_{l_k = 0}^{p_k} {p_1 \choose l_1} \dots {p_k \choose l_k} \theta_L^{p_1 + \dots p_k - (l_1 + \dots + l_k)} (\theta_R - \theta_L)^{l_1 + \dots + l_k} \E\left[U^{l_1}_{\lfloor xN \rfloor+1:N} \dots U^{l_k}_{\lfloor xN \rfloor+k:N}\right].\nonumber
\end{align}
To lighten the notation we write $P_k := p_1 + \dots p_k$ and $L_k = l_1 + \dots + l_k$. We now use the expression in {\color{black}\eqref{eq:exp_products_order_stat}} for the expectation of the standard order statistics.

\begin{align} \label{eq:std_order_stat_loc_equi}
&\E\left[U^{l_1}_{\lfloor xN \rfloor+1:N} \dots U^{l_k}_{\lfloor xN \rfloor+k:N}\right]\\
&\qquad = \frac{N!}{\lfloor xN \rfloor!} \cdot \left( \frac{1}{l_1 + \lfloor xN \rfloor +1} \dots \frac{1}{L_k + \lfloor xN \rfloor +k} \right) \cdot \frac{(L_k + \lfloor xN \rfloor + k)!}{(L_k+N)!}.\nonumber
\end{align}

Notice that 

\be \label{eq:loc_equi_part1}
\frac{N!}{(L_k+N)!} = \frac{1}{(N+1)(N+2) \dots (N + L_k)}
\ee

and 

\begin{align} \label{eq:loc_equi_part2}
&\frac{(L_k+ \lfloor xN \rfloor+k)!}{\lfloor xN \rfloor!} \cdot \left( \frac{1}{l_1 + \lfloor xN \rfloor +1} \dots \frac{1}{L_k + \lfloor xN \rfloor +k} \right)\\
&\qquad = (\lfloor xN \rfloor + 1) \dots (\lfloor xN \rfloor + L_k) \cdot \left(\frac{L_k + \lfloor xN \rfloor + 1 }{l_1 + \lfloor xN \rfloor + 1 } \dots \frac{L_k + \lfloor xN \rfloor + k }{L_k + \lfloor xN \rfloor + k } \right).\nonumber
\end{align}

Combining \eqref{eq:loc_equi_part1} and \eqref{eq:loc_equi_part2} with \eqref{eq:std_order_stat_loc_equi} yields

\begin{align} \label{eq:order_expectation_stdU}
&\E\left[U^{l_1}_{\lfloor xN \rfloor+1:N} \dots U^{l_k}_{\lfloor xN \rfloor+k:N}\right]\\
&\qquad= \frac{\lfloor xN \rfloor + 1}{N + 1} \dots \frac{\lfloor xN \rfloor + L_k}{N + L_k} \cdot \left(\frac{L_k + \lfloor xN \rfloor + 1 }{l_1 + \lfloor xN \rfloor + 1 } \dots \frac{L_k + \lfloor xN \rfloor + k }{L_k + \lfloor xN \rfloor + k } \right)\nonumber\\
&\qquad= \left(x + O(N^{-1})\right)^{L_k} \cdot \left(1 + O(N^{-1})\right)^k = x^{L_k} + O(N^{-1}).\nonumber
\end{align}

Substituting \eqref{eq:order_expectation_stdU} into the expansion in \eqref{eq:loc_equi_expansion} gives

\begin{align}
    &\E \Bigg[ D_N \Bigg(\eta, \sum_{j=1}^k p_j \delta_{\lfloor xN \rfloor + j}\Bigg) \Bigg] \\
    &\qquad = \sum_{l_1 = 0}^{p_1} \dots \sum_{l_k = 0}^{p_k} {p_1 \choose l_1} \dots {p_k \choose l_k} \theta_L^{P_k - L_k} (\theta_R - \theta_L)^{L_k} \left(x^{L_k} + O(N^{-1})\right)\nonumber\\
    &\qquad = (\theta_L + (\theta_R - \theta_L) x )^{P_k} + O(N^{-1}) = \rho_{\theta_L,\theta_R}(x)^{p_1+ \dots + p_k} + O(N^{-1}).
\end{align}
\epr

\section{Law of large numbers for fields of local functions}
In this section we will study fields associated to a general local function, and show their convergence to a deterministic limit.
The density field
is defined as a random measure on $[0,1]$ via
\be
X_N= \frac1N\sum_{i=1}^N \eta_{i,N}\  \delta_{\tfrac{i-1}{N+1}}
\ee 
It acts on smooth test functions
$\phi: [0,1]\to\R$ as follows
\be\label{xnphi}
X_N(\phi)=\frac1N \sum_{i=1}^N \eta_{i,N}\ \phi \left(\frac{i-1}{N+1} \right).
\ee 
For a general local function $g$, we define the field associated to $g$ via its action on smooth test functions:
\be\label{xng}
X_N(g;\phi)= \frac{1}{N} \sum_{i=0}^{N-k} (\tau_i g)[\eta] \phi \left(\frac{i}{N+1} \right).
\ee
Notice that for $g(\eta)=\eta_1$ this reduces to \eqref{xnphi}.

\bt\label{LLNTHM}
Let $\phi: [0,1] \rightarrow \R$ be differentiable with bounded derivative and let $g: \N^N \rightarrow \R$ be local with dependence set $\{1,\dots,k\} \subset \N$, $k\in \N$. \\

\noindent
Assume $g$ is such that there exists a $B(\theta_L,\theta_R) \in \R$ for which the following two conditions hold. For arbitrary $i \in {0, 1, \dots, N-k}$ and $\theta_1, \theta_2, \dots, \theta_N$ $\in [\theta_L,\theta_R]$ such that $\theta_L \leq \theta_1 \leq \theta_2 \leq \dots \leq \theta_{N} <\theta_R $:

\begin{itemize}
    \item [i)] Denote $Y_{i,N} := (\tau_i g)[\eta] - \E[(\tau_i g)[\eta]| \Theta = \theta]$, then, for some $\epsi > 0$,
    \be
    \E \left[ |Y_{i,N}|^{2+\epsi}\right] \leq B(\theta_L,\theta_R) < \infty.
    \ee

    \item[ii)]
    \be
    {\color{black}\left|\E \left[g(\eta)| \Theta = \theta\right]\right|} \leq B(\theta_L,\theta_R) < \infty.
    \ee
\end{itemize}

\noindent
Then the following strong law of large numbers holds,
    \begin{equation}\label{eq:slln_harmonic}
        X_N(g;\phi) = \frac{1}{N} \sum_{i=0}^{N-k} (\tau_i g)[\eta] \phi \left(\frac{i}{N+1} \right) \rightarrow \int_0^1 h(\rho_{\theta_L,\theta_R}(x) ) \phi(x) dx
    \end{equation}
    almost surely as $N \rightarrow \infty$. Here, the function $h: [\theta_L, \theta_R] \rightarrow \R$ is defined as

    \be\label{hgro}
    h(\rho) = h(g, \rho)= \E \left[g(\eta)| \Theta = (\rho, \rho, \ldots, \rho)\right].
    \ee
\et
\br
As will be clear from the proof, we will show that
\[
\frac{1}{N} \sum_{i=0}^{N-k} (\tau_i g)[\eta] \phi \left(\frac{i}{N+1} \right)
\approx 
\frac{1}{N} \sum_{i=0}^{N-k} \E((\tau_i g)[\eta]) \phi \left(\frac{i}{N+1} \right)
\approx
\int_0^1 h(\rho_{\theta_L,\theta_R}(x)) \phi(x) dx
\]
which explains why we call it ``law of large numbers''.
\er 

{\color{black}\bpr
\noindent\textit{Step 1: Decomposition.}
Recall that
\[
Y_{i,N}:=(\tau_i g)[\eta]-\E[(\tau_i g)[\eta]\mid\Theta].
\]
We decompose the sum in \ref{eq:slln_harmonic} into two terms. The first we will refer to as the centered term and the second as the conditional-mean term:
\begin{align}
\frac{1}{N}\sum_{i=0}^{N-k}(\tau_i g)[\eta]\,
\phi\left(\frac{i}{N+1}\right)
={}&\frac{1}{N}\sum_{i=0}^{N-k}Y_{i,N}\,
\phi\left(\frac{i}{N+1}\right)\nonumber\\
&+\frac{1}{N}\sum_{i=0}^{N-k}
\E[(\tau_i g)[\eta]\mid\Theta]\,
\phi\left(\frac{i}{N+1}\right).
\label{eq:lln_decomposition}
\end{align}

\noindent\textit{Step 2: The centered term.}
Fix $r\in\{0,\ldots,k-1\}$ and consider the indices
$i\equiv r\pmod{k}$.  For two distinct indices in this class, the
corresponding functions depend on disjoint sets of the random variables $\eta_1, \eta_2, \ldots, \eta_N$.
Hence, conditionally on $\Theta$, the random variables
\[
Y_{i,N}\phi\left(\frac{i}{N+1}\right),
\qquad 0\leq i\leq N-k,\quad i\equiv r\pmod{k},
\]
are independent and have mean zero.  The conditional Marcinkiewicz--Zygmund
inequality gives
\begin{align}
&\E\left[
\left|
\sum_{\substack{0\leq i\leq N-k\\ i\equiv r\ ({\rm mod}\ k)}}
Y_{i,N}\phi\left(\frac{i}{N+1}\right)
\right|^{2+\epsi}
\Biggm|\Theta
\right]\nonumber\\
&\qquad\leq
C\,\E\left[
\left(
\sum_{\substack{0\leq i\leq N-k\\ i\equiv r\ ({\rm mod}\ k)}}
Y_{i,N}^2\phi\left(\frac{i}{N+1}\right)^2
\right)^{1+\epsi/2}
\Biggm|\Theta
\right]
\leq C N^{1+\epsi/2}.
\label{eq:lln_residue_moment}
\end{align}

For the second inequality we use that
$(\sum_{j=1}^m a_j)^{1+\epsi/2}
\leq m^{\epsi/2}\sum_{j=1}^m a_j^{1+\epsi/2}$ for $a_j\geq0$, and assumption i). The latter is a consequence of Jensen's inequality. Here and below, $C$ is used to denote a constant which is independent of $N$, $r$, and the realization of $\Theta$. 

If the absolute value of the full centred sum is larger than $tN$,
then for at least one residue class the corresponding absolute value
is larger than $tN/k$.  The conditional Markov
inequality therefore yields
\begin{align}
&\P\left(
\left|
\frac{1}{N}\sum_{i=0}^{N-k}
Y_{i,N}\phi\left(\frac{i}{N+1}\right)
\right|>t
\Biggm|\Theta
\right)\nonumber\\
&\qquad\leq
\sum_{r=0}^{k-1}
\P\left(
\left|
\sum_{\substack{0\leq i\leq N-k\\ i\equiv r\ ({\rm mod}\ k)}}
Y_{i,N}\phi\left(\frac{i}{N+1}\right)
\right|>\frac{tN}{k}
\Biggm|\Theta
\right)
\leq C_t N^{-1-\epsi/2}.
\label{eq:lln_summable_bound}
\end{align}
The bound is uniform in $\Theta$. $C_t$ denotes a constant which only depends on $t$.  After taking the expectation with respect to $\Theta$, its
right-hand side is summable in $N$.  The Borel--Cantelli lemma, first
for $t=1/m$, $m\in\N$, then gives
\[
\frac{1}{N}\sum_{i=0}^{N-k}
Y_{i,N}\phi\left(\frac{i}{N+1}\right)
\xrightarrow{a.s.}0.
\]

\noindent\textit{Step 3: The conditional-mean term.}
For fixed value of $\Theta$,
\[
\E[(\tau_i g)[\eta]\mid\Theta]
=
\int g(\xi_1,\ldots,\xi_k)
\prod_{j=1}^k d\nu_{\Theta_{i+j,N}}(\xi_j).
\]
As a function of the $k$ hidden parameters, this expression is local,
bounded and continuous by the assumptions of the theorem and the
geometric form of the measures $\nu_\theta$.  Proposition
\ref{prop:lln_theta} therefore gives
\[
\frac{1}{N}\sum_{i=0}^{N-k}
\E[(\tau_i g)[\eta]\mid\Theta]\,
\phi\left(\frac{i}{N+1}\right)
\xrightarrow{a.s.}
\int_0^1 h(\rho_{\theta_L,\theta_R}(x))\phi(x)\,dx.
\]
Here replacing $i/(N+1)$ by $i/N$ does not affect the limit, since
$\phi$ has bounded derivative and the conditional expectations are
uniformly bounded.  Combining this convergence with
\eqref{eq:lln_decomposition} and Step 2 proves the result.

\epr}

\br
Theorem \ref{LLNTHM} shows that on the level of the law of large numbers, the field of a general local function can be replaced by a function of the density field.
To explain this, for 
$\epsilon N< i< (1-\epsilon) N $ denote by
\[
A_{i,N,\epsilon}(\eta)=\frac{1}{2\epsilon N+ 1}\sum_{|j-i|<N\epsilon} \eta_j
\]
the empirical density in the macroscopic small block
$[i-\epsi N, i+\epsi N]$.
 
Then Theorem \ref{LLNTHM} shows that we can make the following replacement in the limit $N\to\infty, \epsi\to 0$.
\begin{align}\label{bawafa}
X_N(g;\phi) &= \frac{1}{N} \sum_{i=1}^{N-k} (\tau_i g)[\eta] \phi \left(\frac{i}{N+1} \right)
\nonumber\\
& \approx
\frac{1}{N}
\sum_{i=N\epsilon}^{N-N\epsilon} h(A_{i,N,\epsilon}(\eta))
\phi \left(\frac{i}{N+1} \right)
\nonumber\\
& \approx
\int_{\epsilon}^{1-\epsilon} h(X_N*i_{\epsilon}(x)) \phi(x) dx
\end{align}
Here, $i_\epsi= {\color{black} (2 \epsi)^{-1} 1_{[-\epsilon,\epsilon]}}$ and $*$ denotes convolution, i.e., for $\epsilon\leq x\leq 1-\epsilon$:
\[
\mu*i_\epsilon (x)= \int_0^1 \mu(dy) i_{\epsi} (x{\color{black}-}y),
\]
and $\approx$ means that the difference between the two quantities tends to zero almost surely when $N\to\infty$, followed by $\epsi\to 0$.
\er
\br
Because the law of large numbers in Theorem \ref{LLNTHM}
holds for all test functions $\phi$, we can reformulate it
as an almost sure weak convergence of profiles.
More precisely, for a local function $g$ as in the statement of Theorem \ref{LLNTHM}, we define
its associated empirical profile via
\be\label{infdimlln}
\Xi(g,N;\eta) =\frac{1}{N} \sum_{i=0}^{N-k} (\tau_i g)[\eta] \delta_{\frac{i}{N+1} }.
\ee
Then Theorem \ref{LLNTHM} can be reformulated as 
follows. When $N\to\infty$, the random measures
$\Xi(g,N;\eta)$ almost surely weakly converge to the deterministic measure given by 
$h(\rho_{\theta_L,\theta_R}(x)) dx$.
Generalizing even more, we can also vary the function $g$ and consider the 
random distribution which maps the pair $(g, \phi)$ to
\[
\frac{1}{N} \sum_{i=0}^{N-k} (\tau_i g)[\eta] \phi\left(\frac{i}{N+1} \right)
\]
Then we obtain that this random distribution almost surely weakly converge to the deterministic measure given by 
$h(\rho_{\theta_L,\theta_R}(x)) dx$, where $h$ is given by \eqref{hgro}.
\er 

\section{Central limit theorem for fluctuation fields of local functions}

In this section we look at the fluctuation fields, i.e., the fluctuations around the non-equilibrium steady state profile on central limit scale.
In equilibrium, the stationary state is reversible and is a product measure, and as a consequence these fluctuation fields converge to white noise.
In non-equilibrium (i.e., different reservoir parameters
$\theta_L\not=\theta_R$), there is a white noise contribution (corresponding to the local equilibrium approximation of the stationary state), and a contribution coming from the fluctuation of the random parameters, which are responsible for the long-range correlations in the NESS.
As we will see below, the latter behaves as a Brownian bridge (cf.\ also \eqref{prop:clt_Theta}).

\bt
Let $\phi:[0,1] \rightarrow \mathbb{R}$ be a smooth function and let $g: \N^N\rightarrow \R$ be a local function such that 

\be
\E\left[\big((\tau_n g)[\eta] - \E((\tau_n g)[\eta]| \Theta = \theta)\big)^4\right] < \infty
\ee

and

\be
\cov\Big( (\tau_n g)[\eta] - \E((\tau_n g)[\eta]| \Theta = \theta) \,\, , \,  (\tau_m g)[\eta] - \E((\tau_m g)[\eta]| \Theta = \theta)\Big) < \infty
\ee

for all $n,m \in \{0, 1, \ldots,N-k\}$ and $\theta \in [\theta_L, \theta_R]^N$. Then we have
    \begin{equation}\label{mainclt}
        \lim_{N \rightarrow \infty} \frac{1}{\sqrt{N}} \sum_{n=0}^{N-k} \big((\tau_n g)[\eta] - \E((\tau_n g)[\eta]) \big) \phi \Big(\frac{n}{N+1} \Big) \stackrel{\text{d}}{=} \mathcal{N}(0,\sigma_T^2 + \sigma_E^2).
    \end{equation}
Here $\sigma_T^2$ and $\sigma_E^2$ are given by 

    \begin{align}\label{varii}
    &\sigma_T^2 = (\theta_R-\theta_L)^2 \int_{0}^{1} \int_0^1 [s \wedge t -st] \phi(s) \phi(t) h'(\rho_{\theta_L, \theta_R}(s)) h'(\rho_{\theta_L, \theta_R}(t)) ds dt\\
    &\sigma_E^2 = \int_0^1 V(\rho_{\theta_L, \theta_R}(x)) \phi^2(x) dx.
    \end{align}

Here, the functions $h:[\theta_L, \theta_R] \rightarrow \R$ and $V:[\theta_L, \theta_R] \rightarrow \R$ are defined as
    
    \begin{align}
    &h(\theta) = \E[g(\eta_1, \dots, \eta_k)|\Theta_{1,N} = \Theta_{2,N} = \dots = \Theta_{k,N} = \theta],\\
    &\nonumber\\
    &V(\theta) = \sum_{m=1}^{2k-1} \cov(g(\eta_{k}, \dots, \eta_{2k-1}), g(\eta_{m}, \dots, \eta_{m+k-1})|\Theta_{1,N} = \Theta_{2,N} = \dots = \Theta_{k,N} = \theta)
    \end{align}

In particular, when $g(\eta)=\eta_1$ we obtain the following result for the density fluctuation field.
\be\label{bridge}
\lim_{N \rightarrow \infty} \frac{1}{\sqrt{N}} \sum_{n=1}^{N} \big(\eta_n - \E(\eta_n) \big) \phi \Big(\frac{n}{N+1} \Big) = (\theta_R-\theta_L)\int_0^1 \phi(s) \overline{B}(s) ds + \mathcal{W}(\phi)
\ee
Here the limit is in distribution,
$\{\overline{B}(s): 0\leq s\leq 1\}$ denotes a Brownian bridge, and $\mathcal{W}(\phi)$ denotes an independent white noise with variance
\[
\sigma_E^2= \int_0^1 \var(\eta_1|\Theta_{1,N}=\rho_{\theta_L, \theta_R}(x))\ \phi^2(x) dx.
\]
\et
\br
The term $\mathcal{N} (0, \sigma_E^2)$ corresponds to the variance of the fluctuation field in the local equilibrium product measure
\be\label{locprod}
\otimes_{n=1}^N \nu_{\rho_{\theta_L,\theta_R}(n/(N+1))}
\ee
where we recall the notation $\rho_N(x)=\theta_L + x(\theta_R-\theta_L)$, $x\in [0,1]$.
This white noise term is also present in equilibrium, i.e., when $\theta_R=\theta_L$.
The term $\mathcal{N} (0, \sigma_T^2)$ corresponds to the contribution of the (long-range) correlations of the random parameters $\Theta_{n,N}$, $n=1, \ldots, N$. This term vanishes in equilibrium, i.e., when $\theta_L=\theta_R$.
The result for the fluctuation field of a general local function can then be restated as follows.
\begin{align}\label{bridge_2}
&\lim_{N \rightarrow \infty} \frac{1}{\sqrt{N}} \sum_{n=0}^{N-k} \big((\tau_n g)[\eta] - \E((\tau_n g)[\eta]) \big) \phi \Big(\frac{n}{N+1} \Big)\\
& \hspace{3cm} = (\theta_R-\theta_L)\int_0^1 
h'(\rho_{\theta_L, \theta_R}(s)) {\color{black} \phi(s)} \overline{B}(s) ds + \mathcal{W}_g(\phi).\nonumber
\end{align}
Here the white noise contribution $\mathcal{W}_g(\phi)$ corresponds to the limiting variance of the fluctuation field in the local equilibrium product measure \eqref{locprod} (term $\si_E^2$ in \eqref{varii}), whereas the Brownian bridge contribution corresponds to the correlations of the random parameters. We observe that this Brownian bridge contribution differs from the corresponding term in \eqref{bridge} only by the addition of the factor 
$h'(\rho_{\theta_L, \theta_R}(s)) $, which is a function of the local density field.
\er 

\bpr
To prove \eqref{bridge} given \eqref{mainclt} , notice that if
$g(\eta)=\eta_1$, $h(x)=x$, and therefore
\[
\sigma_T^2 = (\theta_R-\theta_L)^2 \int_{0}^{1} \int_0^1 [s \wedge t -st] \phi(s) \phi(t)  ds dt
\]
which is precisely the variance of the random variable 
$(\theta_R-\theta_L)\int_0^1 \phi(s) \overline{B}(s) ds$, whereas the white noise term $\mathcal{W}(\phi)$ produces the $\sigma_E^2$ term in the variance.

Next we prove \eqref{mainclt} via the convergence of the characteristic function.

\begin{align} \label{eq:begin_CLT}
\E \Bigg[&\exp \left\{ \frac{i t}{\sqrt{N}} \sum_{n=0}^{N-k} \big((\tau_n g)[\eta] - \E((\tau_n g)[\eta]) \big) \phi \Big(\frac{n}{N+1} \Big) \right\} \Bigg] =\\
\E \Bigg[&\exp \left\{ \frac{i t}{\sqrt{N}} \sum_{n=0}^{N-k} \big(\E[(\tau_n g)[\eta]|\Theta] - \E((\tau_n g)[\eta]) \big) \phi \Big(\frac{n}{N+1} \Big) \right\}\nonumber\\
&\cdot \E\left( \exp \left\{ \frac{i t}{\sqrt{N}} \sum_{n=0}^{N-k} \big((\tau_n g)[\eta] - \E((\tau_n g)[\eta]|\Theta) \big) \phi \Big(\frac{n}{N+1} \Big) \right\} \Bigg| \Theta \right) \Bigg]. \nonumber
\end{align}

It turns out that the expectation conditional on $\Theta$
in the rhs of \eqref{eq:begin_CLT} will converge almost surely (in the realization of $\Theta$) to the characteristic function of a normal random variable with variance $\sigma_E^2$. To show this, we use a version of the central limit theorem which holds for a triangular scheme of weakly dependent random variables \cite{neumann2013central}. Let $X_{N,n}$ denote 

\be
X_{N,n} :=  \big((\tau_n g)[\eta] - \E((\tau_n g)[\eta]|\Theta) \big) \phi \Big(\frac{n}{N+1} \Big)
\ee

The central limit theorem in \cite{neumann2013central} states that 

\be 
\frac{1}{\sqrt{N}}(X_{N,1} + \dots + X_{N,N-k}) \xrightarrow{d} \mathcal{N}(0, \sigma^2),
\ee

with $\sigma^2$ as in \eqref{eq:sigma_Neumann_CLT}, if the following four conditions are satisfied:
\begin{itemize}
    \item[1)] $\E[ X_{N,m} |\Theta] = 0$ and $\frac{1}{N}\sum_{n=1}^{N-k} \E[X_{N,n}^{2}|\Theta] \leq v_0$, for all $N, m$ and some $v_0 < \infty$.
    
    \item[2)] The variance of the sum converges:
    
    \be\label{eq:sigma_Neumann_CLT}
    \sigma_N^2 =\frac{1}{N} \text{var}(X_{N,1} + \dots + X_{N,N}| \Theta) \xrightarrow{N \to \infty} \sigma^2 \in [0,\infty].
    \ee
    
    \item[3)] The Lindeberg condition  holds: for all $\epsilon>0$,
    
    \be \label{eq:lindeberg}
    \sum_{n=1}^{N-k} \E \left[\frac{1}{N}X_{N,n}^{2} \mathds{1} \left(|X_{N,n}| > \epsilon \sqrt{N}\right)\Big| \Theta \right] \xrightarrow{N \to \infty} 0.
    \ee
    
    \item[4)] There exists a summable sequence $(B_r)_{r \in \mathbb{N}}$ such that, for all $u \in \mathbb{N}$ and all indices $1 \leq s_1 < s_2 < \dots < s_u < s_u + r = t_1 \leq t_2 \leq N-k$, the following upper bounds for covariances hold true: for all measurable functions $f : \mathbb{R}^{u} \to \mathbb{R}$ with $\| f \|_{\infty} = \sup_{x \in \mathbb{R}^{u}} |f(x)| \leq 1$,
    \begin{align} \label{eq:mixing1}
    \big| \cov (f(X_{N,s_1}, \dots, X_{N,s_u}) X_{N,s_u}, &X_{N,t_1} | \Theta) \big|\\
    & \leq \left(\E[ X_{N,s_u}^{2}|\Theta] + \E[ X_{N,t_1}^{2}|\Theta] + \frac{1}{N-k} \right) B_r \nonumber
    \end{align}
    and
    \begin{align} \label{eq:mixing2}
    \big| \cov (f(X_{N,s_1}, \dots, X_{N,s_u}), X_{N,t_1} &X_{N,t_2} |\Theta) \big|\\
    &\leq \left(\E[ X_{N,t_1}^{2}| \Theta] + \E[ X_{N,t_2}^{2}|\Theta] +  \frac{1}{N-k} \right) B_r.\nonumber
    \end{align}
\end{itemize}

We verify these conditions. The first condition is satisfied because $\E[X_{N,n}^{2}|\Theta]$ is bounded, as it is a continuous function of $\Theta_{n-1},...,\Theta_{k+n-1}$ on the closed interval $[\theta_L,\theta_R]^k$.\\

For the second condition we notice that 
\begin{align}
    &\frac{1}{N} \var \left(\sum_{n=1}^{N-k} X_{N,n}\Big| \Theta \right) = \frac{1}{N} \sum_{n,m=1}^{N-k} \cov(X_{N,n},X_{N,m}|\Theta)\\
    &\qquad=  \frac{1}{N} \sum_{n=1}^{N-k} \quad  \sum_{m= (n+1-k) {\color{black} \vee} 1 }^{(n-1+k) {\color{black}\wedge} (N-k)} \cov(X_{N,n} (\eta_n, \dots, \eta_{n+k-1}),X_{N,m} (\eta_m, \dots, \eta_{m+k-1})|\Theta) \nonumber\\
    &\qquad =\frac{1}{N} \sum_{n=k}^{N+1-2k} \quad  \sum_{m= n+1-k}^{n-1+k} \cov(X_{N,n} (\eta_n, \dots, \eta_{n+k-1}),X_{N,m} (\eta_m, \dots, \eta_{m+k-1})|\Theta) + {\color{black} o(1)} \nonumber\\
    &\qquad= \frac{1}{N} \sum_{n=k}^{N+1-2k} \phi \left(\frac{n}{N-k}\right)^2  \hat{V}(\Theta_{n-k+1,N}, \dots, \Theta_{n+2k-1,N}) + {\color{black} o(1)} ,\nonumber
\end{align}

with $\hat{V}:[\theta_L,\theta_R]^{2k-1} \rightarrow \R$ defined as 

\begin{align}
    &\hat{V}(\theta_1, \dots, \theta_{2k-1} )\\
    &\qquad= \sum_{m=1}^{2k-1} \cov(g(\eta_{k}, \dots, \eta_{2k-1}), g(\eta_{m}, \dots, \eta_{m+k-1})| \Theta_{1,N} = \theta_1, \dots, \Theta_{2k-1,N} = \theta_{2k-1}).\nonumber
\end{align}

Using Proposition \ref{prop:lln_theta}, we find

\begin{align}
&\frac{1}{N} \var \left(\sum_{n=1}^{N-k} X_{N,n} \Big| \Theta \right)\\
&\qquad= \frac{1}{N} \sum_{n=k}^{N+1-2k} \phi \left(\frac{n}{N-k}\right)^2  \hat{V}(\Theta_{n-k+1,N}, \dots, \Theta_{n+2k-1,N}) + {\color{black}o}(1)\nonumber\\
&\qquad \rightarrow \int_0^1  V(\rho_{\theta_L,\theta_R}(x)) \phi (x)^2 dx = \sigma^2_E \nonumber
\end{align}

as $N \rightarrow \infty$. \\

For the third condition, we bound each term

\begin{align}
    \E \left[X_{N,n}^{2} \mathds{1} \left(|X_{N,n}| > \epsilon \sqrt{N}\right) \Big| \Theta \right] &\leq \E\left[X_{N,n}^{4} \Big| \Theta \right]^{\frac{1}{2}} \E\left[ \mathds{1}^2 \left(|X_{N,n}| > \epsilon \sqrt{N}\right) \Big| \Theta \right]^{\frac{1}{2}}\\
    &= \E\left[X_{N,n}^{4} \Big| \Theta \right]^{\frac{1}{2}} \Pr \left( |X_{N,n}| > \epsilon \sqrt{N} \Big| \Theta \right)^{\frac{1}{2}}\nonumber\\
    &\leq \frac{\E\left[X_{N,n}^{4} \Big| \Theta \right]}{\epsilon^2 N}. \nonumber
\end{align}

Notice that $\E [X_{N,n}^{4} | \Theta ]$ is bounded because it is a continuous function of $\Theta_{n,N}, \dots,$ $\Theta_{n+k-1,N}$ on $[\theta_L,\theta_R]^{k}$. We see that the left hand side of \eqref{eq:lindeberg} indeed converges to $0$ as $N\rightarrow \infty$.\\

We show condition four. For readability we suppress the dependence of\\ $f(X_{N,s_1}, \dots, X_{N,s_u})$ on $X_{N,s_1}, \dots, X_{N,s_u}$ and write simply $f$. We find an upper bound for the left hand side of \eqref{eq:mixing1}. Notice that

\begin{align} \label{eq:mixing1_pf}
     \big| \cov (f \cdot X_{N,s_u}, X_{N,t_1} | \Theta) \big| &= \left| \E \left[(f \cdot X_{N,s_u} - \E[f \cdot X_{N,s_u}|\Theta] )X_{N,t_1}\Big|\Theta\right] \right|\\
     &\leq \left| \E \left[f \cdot X_{N,s_u} X_{N,t_1}\Big|\Theta\right] \right| + \Big|\E\left[f \cdot X_{N,s_u}\Big|\Theta\right] \E\left[X_{N,t_1}\Big|\Theta\right]\Big| \nonumber\\
     &\leq \E \left[|f|  |X_{N,s_u}| |X_{N,t_1}|\Big|\Theta\right]\nonumber\\
     &\leq \E \left[|X_{N,s_u}| |X_{N,t_1}|\Big|\Theta \right]\nonumber\\
     &\leq \E\left[|X_{N,s_u}|^2 \Big|\Theta\right]^{\frac{1}{2}} \E\left[|X_{N,t_1}|^2\Big|\Theta\right]^{\frac{1}{2}} \nonumber\\
     &\leq \frac{1}{2} \left(\E\left[|X_{N,s_u}|^2 \Big|\Theta\right] + \E\left[|X_{N,t_1}|^2\Big|\Theta\right]\right).\nonumber
\end{align}

For the third line of \eqref{eq:mixing1_pf} we use that $\E[X_{N,t_1}|\Theta]$ is zero. Then we use that $\norm{f}_\infty = 1$ for the fourth line and Cauchy-Schwarz for the fifth line. Finally we use the fact that $0 \leq (a - b)^2 = a^2 +b^2 - 2ab$ for $a,b \in \R$ to obtain the last line.

Next we also find an upper bound for the left hand side of \eqref{eq:mixing2}.

\begin{align} \label{eq:mixing2_pf}
     \big| \cov (f , X_{N,s_u} X_{N,t_1} | \Theta) \big| &= \left| \E \left[(f  - \E[f |\Theta] ) \cdot (X_{N,s_u} X_{N,t_1} - \E[X_{N,s_u} X_{N,t_1} |\Theta]) \Big|\Theta\right] \right|\\
     &\leq \E \left[|f  - \E[f |\Theta] | \cdot |X_{N,s_u} X_{N,t_1} - \E[X_{N,s_u} X_{N,t_1} |\Theta]| \Big|\Theta\right] \nonumber\\
     &\leq 2\E \left[ |X_{N,s_u} X_{N,t_1} - \E[X_{N,s_u} X_{N,t_1} |\Theta]| \Big|\Theta\right]\nonumber\\
     &\leq 4 \E\left[|X_{N,s_u}|| X_{N,t_1} | \Big| \Theta \right] \nonumber\\
     &\leq 4 \E\left[|X_{N,s_u}|^2 \Big|\Theta\right]^{\frac{1}{2}} \E\left[|X_{N,t_1}|^2\Big|\Theta\right]^{\frac{1}{2}} \nonumber\\
     &\leq 2 \left(\E\left[|X_{N,s_u}|^2 \Big|\Theta\right] + \E\left[|X_{N,t_1}|^2\Big|\Theta\right]\right) \nonumber
\end{align}

In the third line of \eqref{eq:mixing2_pf} we use that $\norm{f}_\infty = 1$. For the fifth line we use Cauchy -Schwarz. In the last line we use again that $a^2+b^2 \geq 2ab$ for $a,b \in \R$.

Now notice that for $r \geq k$, we have

\be
\cov (f \cdot X_{N,s_u}, X_{N,t_1} | \Theta) =  \cov (f , X_{N,s_u} X_{N,t_1} | \Theta) = 0,
\ee

because $X_{N,s_1}, \dots X_{N,s_u}$ are a function of $\eta_1, \dots, \eta_{s_u + k-1}$ while $X_{N,t_1}, X_{N,t_2}$ are a function of $\eta_{s_u + r} , \dots \eta_{t_2 + k-1}$. Therefore, we can define $B_r$ as

\be
B_r := 2 \cdot \mathds{1}(r < k).
\ee

Since $B_r$ is summable in $r$, we have shown the fourth condition. We conclude that   

\be
\frac{1}{\sqrt{N}} \sum_{n=0}^{N-k} \big((\tau_n g)[\eta] - \E((\tau_n g)[\eta]|\Theta) \big) \phi \Big(\frac{n}{N+1} \Big) \xrightarrow{d} \mathcal{N}(0, \sigma_E^2)
\ee

with respect to $\Pr(\cdot|\Theta)$.{\color{black} Hence the corresponding characteristic function converges to $\exp\{-\sigma_E^2 t^2/2\}$,}

\be
\E\left( \exp \left\{ \frac{{\color{black}i} t}{\sqrt{N}} \sum_{n=0}^{N-k} \big((\tau_n g)[\eta] - \E((\tau_n g)[\eta]|\Theta) \big) \phi \Big(\frac{n}{N+1} \Big) \right\} \Bigg| \Theta \right) \rightarrow \exp\left\{\frac{{\color{black}-}\sigma^2_E t^2}{2}\right\},
\ee

{\color{black}almost surely in $\Theta$. As a consequence, \eqref{eq:begin_CLT} gives } 

\begin{align}
\E \Bigg[&\exp \left\{ \frac{{\color{black}i}t}{\sqrt{N}} \sum_{n=0}^{N-k} \big((\tau_n g)[\eta] - \E((\tau_n g)[\eta]) \big) \phi \Big(\frac{n}{N+1} \Big) \right\} \Bigg] =\\
\E \Bigg[&\exp \left\{ \frac{{\color{black}i}t}{\sqrt{N}} \sum_{n=0}^{N-k} \big(\E[(\tau_n g)[\eta]|\Theta] - \E((\tau_n g)[\eta]) \big) \phi \Big(\frac{n}{N+1} \Big) \right\} \exp\left\{\frac{{\color{black}-}\sigma^2_E }{2} t^2\right\} \Bigg]+{o(1)}. \nonumber
\end{align}
Here $o(1)$ denotes a term converging to zero as $N\to \infty$.
Using Proposition \ref{prop:clt_Theta} we find 

\begin{align}
\E \Bigg[&\exp \left\{ \frac{{\color{black}i} t}{\sqrt{N}} \sum_{n=0}^{N-k} \big((\tau_n g)[\eta] - \E((\tau_n g)[\eta]) \big) \phi \Big(\frac{n}{N+1} \Big) \right\} \Bigg] \rightarrow \exp\left\{{\color{black}-} \frac{\sigma^2_T + \sigma^2_E}{2} t^2\right\}.
\end{align}

This concludes the proof.
\epr

\section{Large deviations for fields of local functions }
In \cite{carinci2025large} we proved a large deviation principle for the density field in the NESS of the generalized harmonic models, with a rate function closely related to the rate function obtained in \cite{bertini2007stochastic}, \cite{bertini2}. In this section we generalize this large deviation result to the field of a local function. 
The idea is that under a product measure of the form
$\otimes_{i=1}^N \nu_{\Theta(i/N)}$,
the field of a local function
(cf. \eqref{xng}) satisfies a large deviation principle with a large deviation free energy only depending on the profile of the $\theta$ parameters.
Therefore, because the profile of the $\theta$ parameters satisfies a large deviation principle,  using Varadhan's lemma, one can thus
obtain the large deviation free energy and the large deviation rate function for
the field of a general local function.

\subsection{Large deviations for product measure}
Let $g= g(\eta_1,\ldots, \eta_k)$ denote a bounded local function.
Some of the statements below can be extended to
unbounded functions such as $g(\eta)= \eta_1$, but for the sake of simplicity we restrict to bounded functions, which ensures the existence of moment generating functions such as in \eqref{homfree} below (for all $\lambda\in\R$). 
With a slight abuse of notation
we denote by $\nu_\theta$ the homogeneous product measure
with geometric marginals with parameter $\theta$ (cf. \eqref{geodi}). Under this measure, 
the random variables 
\be
A(N,g)=\frac{1}{N}\sum_{i=1}^{N-k} \tau_i g 
\ee 
satisfy a large deviation principle with rate function
\be\label{homrate}
I(\theta,x; g)= \sup_{\lambda\in\R} (\lambda x - F(\theta, \lambda;g))
\ee
where, in turn,  the large deviation free energy $F(\theta,\lambda;g)$
is given by
\be\label{homfree}
F(\theta,\lambda; g)= \lim_{N\to\infty} 
\frac1N\log\E_{\nu_\theta} e^{\lambda\sum_{i=1}^{N-k} \tau_i g }.
\ee
The existence of the limit \eqref{homfree} as well as the large deviation principle with rate function $I(\theta, \cdot; g)$
is an application of standard Gibbs theory, see e.g.
\cite{georgii}.

The following proposition generalizes this statement for 
inhomogeneous product measures with a given profile of parameters.
\color{black}
More precisely, consider
a non-decreasing c\`{a}dl\`{a}g function $\Theta: [0,1]\to \R^+$ and
the associated product measure
$ \nu_\Theta= \otimes_{i=1}^N \nu_{\Theta(i/(N+1))}$.
\color{black}
Consider then $\phi: [0,1]\to\R$ another continuous
function, and consider
\be\label{inhomav}
A(N,g, \phi)=\frac{1}{N}\sum_{i=1}^{N-k} (\tau_i g) \phi \left(\tfrac{i}{N+1}\right)
\ee
Then we have the following.

\bp\label{waki}
Under the measure $\nu_\Theta= \otimes_{i=1}^N \nu_{\theta(i/(N+1))}$, we have that 
\be\label{inhomfreeen}
\caF(\Theta,\phi;g):=\lim_{N\to\infty}\frac1N\log 
\E_{\nu_\Theta}(e^{N A(N,g, \phi)})
=
\int_0^1 F(\Theta(x), \phi(x); g) dx,
\ee
where $F(\Theta(x), \phi(x); g)$ is given by
\eqref{homfree}.
As a consequence, the random profile
\be\label{profibol}
\Xi(N,g)= \frac{1}{N}\sum_{i=1}^{N-k} \tau_i g \delta_{i/(N+1)}
\ee
satisfies the large deviation principle with rate function given by
\be\label{inhomrate}
\caI(\Theta, \mu(x) dx; g)= 
\int_0^1 I(\Theta (x),  \mu(x); g)\  dx
\ee
and $\caI(\Theta, \mu; g)= \infty$ whenever $\mu$ is not absolutely continuous w.r.t. the Lebesgue measure.
\ep 
\bpr
We will give the proof for $\Theta: [0,1] \to\R$ continuous, and at the end indicate how to extend to the c\`{a}dl\`{a}g case.
It suffices to show \eqref{inhomfreeen}, the rate function
\eqref{inhomrate} then follows by Legendre transform.
\color{black}
Indeed, notice that by the boundedness of $g$, the rate function
is finite only for absolutely continuous measures $\mu(x) dx$, because the corresponding free energy function \eqref{inhomfreeen} is finite for all
continuous $\phi: [0,1]\to\R$.
\color{black}
Notice that this would be different if
we considered unbounded $g$ such as e.g. $g(\eta)=\eta_1$, because then $\phi(x)$ has to be less
than $1/\Theta(x)$, otherwise $F(\Theta(x), \phi(x); g)=\infty$.






Fix $\epsi>0$. We choose $u_{1,N}, u_{2,N}, \dots, u_{l,N}\in [0,1]$, $l\in \N$, such that for all $i \in \{1, \dots, l-1\}$, we have

\be\label{dijo}
k < (\lfloor u_{i+1,N} N \rfloor - \lfloor u_{i,N} N \rfloor ) - 2 \epsi N < k + 2
\ee

and 

\be
\lfloor u_{l,N} N \rfloor  -\epsi N \leq N \leq \lfloor u_{l,N} N \rfloor  +\epsi N.
\ee

We denote the $(\epsi N)-$balls around the points $u_{1,N}, \dots u_{l,N}$ by

\be
B_{i,N} := (\lfloor u_i N \rfloor - \epsi N , \lfloor u_i N \rfloor + \epsi N) \cap \{1, \dots,N\}.
\ee
\br\label{epsrem}
Notice that two distinct $\epsilon N$-balls are at distance larger than $k$ and smaller than $k+2$.
This implies that the union of the $l$ balls
covers all the points of $\{1, \ldots, N\}$ except for at most of the order $k/\epsilon$ points.
\er 
We denote the maximal element of each ball by $M_{i,N}$ and the minimal element by $m_{i,N}$. 
Then we have the following chain of equalities
\begin{align}\label{bopio}
    &\frac{1}{N} \log \int \otimes_{i=1}^N \nu_{\Theta(\frac{i}{N})}(d\eta_{i,N}) \exp \left\{\sum_{x=0}^{N-k} g(\eta_{x+1,N}, \dots, \eta_{x+k,N}) \phi\left( \frac{x}{N}\right) \right\}\nonumber\\
    &= \frac{1}{N} \log \int \otimes_{i=1}^N \nu_{\Theta(\frac{i}{N})}(d\eta_{i,N}) \exp \left\{\sum_{j=1}^l \sum_{x \in B_{j,N}} g(\eta_{x+1,N}, \dots, \eta_{x+k,N}) \phi\left( \frac{x}{N}\right) + O(k/\epsilon)\right\} \nonumber\\
    &= \sum_{j=1}^l \frac{1}{N} \log \int \nu_{\Theta ( (m_{j,N})/N )}(d\eta_{m_{j,N},N}) \dots  \nu_{\theta ( (M_{j,N}+k-1)/N )}(d\eta_{M_{j,N}+k-1,N})\nonumber \\
    & \hspace{2.60cm}\exp \left\{ \sum_{x \in B_{j,N}} g(\eta_{x+1,N}, \dots, \eta_{x+k,N}) \phi\left( \frac{x}{N}\right)\right\} + O\left(\tfrac{k}{\epsilon N}\right)\nonumber\\
    &= \sum_{j=1}^l \frac{1}{N} \log \int \nu_{\Theta ( u_j )}(d\eta_{m_{j,N},N}) \dots  \nu_{\Theta ( u_j )}(d\eta_{M_{j,N}+k-1,N})\nonumber \\
    & \hspace{2.60cm}\exp \left\{ \sum_{x \in B_{j,N}} g(\eta_{x+1,N}, \dots, \eta_{x+k,N}) \phi\left( u_j \right)\right\} + o_{\epsi}(1) + 
    O\left(\tfrac{k}{\epsilon N}\right)\nonumber\\
    &= \sum_{j=1}^l (2 \epsi) F(\Theta(u_j), \phi(u_j);g) + o_{\epsi}(1) + 
    O\left(\tfrac{k}{\epsilon N}\right).
\end{align}
Here in the first step we used remark \ref{epsrem}, combined with the boundedness of $g$.
In the second step we used that the sums over $B_{j,N}$ are
independent by \eqref{dijo}.
In the third step we used that the measure $\nu_\theta$ depends continuously on the parameter $\theta$; here $o_\epsilon(1)$ denotes a term which converges to zero as $\epsilon\to 0$.

By taking the limit $N \rightarrow \infty$ followed by the limit $\epsi \rightarrow 0$ we obtain \eqref{inhomfreeen} because the sum in the rhs of \eqref{bopio} is a Riemann sum approximation of the integral
\eqref{inhomfreeen}.
\color{black}
To extend to the c\`{a}dl\`{a}g case, notice that $\Theta: [0,1]\to\R$ can have at most a countable number of discontinuities, hence we can still choose
the points $u_{1,N}, u_{2,N}, \dots, u_{l,N}\in [0,1]$ to be continuity points of $\Theta$.
\color{black}

\epr
\subsection{Large deviations for mixed product states}
Notice that the precise form $\otimes_{i=1}^N \nu_{\Theta(i/N)}$ of the product measure in Proposition \ref{waki} can be relaxed.
Indeed if the product measure is of the form
$\otimes_{i=1}^N \nu_{\theta_N(i)}$ where $\theta_N(i)$ is such that its associated sample path
\[
\Theta_N(t):= \theta_N(\lfloor tN\rfloor)
\]
converges to a limit $\Theta(t)$ in the Skorokhod topology (of c\`{a}dl\`{a}g paths $\gamma: [0,1]\to \R^+$) then the result of Proposition \ref{waki} still holds, with identical proof, when the expectation over
$\nu_\Theta$ is replaced by the expectation over
$\otimes_{i=1}^N \nu_{\theta_N(i)}$.

We can then combine  Proposition \ref{waki} and Proposition \ref{prop:LDP_order_stat} with Varadhan's lemma to deal with mixed product states, provided the sample path of the random parameters satisfies a large deviation principle. 
\color{black}
Indeed, notice that integral $\int_0^1 F(\Theta(x), \phi(x); g) dx$,
appearing in the rhs of \eqref{inhomfreeen}, is a continuous function of $\Theta$ in the Skorokhod topology.
\color{black}
As a consequence, we have the following large deviation result. Its proof, given the result of 
Proposition \ref{waki}, is completely analogous to the proof given in \cite{carinci2025large} for the large deviations of the density profile, and it is therefore left to the reader.

\bc \label{cor:var_expression}
Let $g:\N^k \rightarrow \R$ be a local and bounded function.
Define $\Theta_N(i/N)= \Theta_{i,N}$ where $\Theta_{i,N}$ is defined in \eqref{ordthet}, i.e., the $i$-th order statistic of $N$ independent uniforms on $[\theta_L,\theta_R]$.
We then consider the corresponding mixed 
product state $\mu_N$ as in \eqref{mixprod2}.

Then we have the following results.
\ben 
\item Large deviation free energy for mixed product measure.
\begin{align}
&\lim_{N \rightarrow \infty} \frac{1}{N} \log \E \int \otimes_{i=1}^N \nu_{\Theta_N(\frac{i}{N})}(d\eta_{i,N}) \exp (N A(N,g,\phi))\\
&\qquad = \sup_{\theta \in A_{0,1}(\theta_L,\theta_R)}\left[ \int_0^1 F(\theta(x), \phi(x); g) dx - J(\theta) \right] \nonumber
\end{align}
Here $J(\theta)$ is the rate function for sample paths of uniform order statistics, defined in 
\eqref{bonki}, and $A_{0,1}(\theta_L,\theta_R)$
denotes the closed set of non-decreasing c\`{a}dl\`{a}g functions $f : [0, 1] \rightarrow \R$ such that $f(x) \geq \theta_L$ and $f(1) = \theta_R$. 

\item Large deviation principle for the field of a local function under the mixed product measure.
Under $\mu_N$, the empirical profile of the local function $\Xi(N,g)$ (defined in \eqref{profibol}) satisfies a large deviation principle with rate function
\[
\caI (\mu(x) dx)=
\inf_{\Theta\in A_{0,1} (\theta_L,\theta_R)}
\left(\caI(\Theta, \mu(x) dx; g) + J(\Theta)\right)
\]
\een
\ec
\subsection*{Acknowledgement}
Berend van Tol was supported by the Peter Paul Peterich
Foundation via the TU Delft University Fund.
We thank Cristian Giardin\`{a}
for several inspiring discussions.\\

\noindent
\textbf{Data Availability:}  Data sharing not applicable to this article as no datasets were generated or analysed during the current study.\\
\textbf{Declarations Competing interests:} The authors have no relevant financial or non-financial interests to disclose.

\bibliographystyle{plain}
\bibliography{bibliography}

@book{georgii,
  title={Gibbs measures and phase transitions},
  author={Georgii, Hans-Otto},
  volume={9},
  year={2011},
  publisher={Walter de Gruyter}
}

@article{bertini2,
  title={Large deviations for a stochastic model of heat flow},
  author={Bertini, Lorenzo and Gabrielli, Davide and Lebowitz, Joel L},
  journal={Journal of statistical physics},
  volume={121},
  number={5},
  pages={843--885},
  year={2005},
  publisher={Springer}
}

@book{shorack2009empirical,
  title={Empirical processes with applications to statistics},
  author={Shorack, Galen R. and Wellner, Jon A.},
  year={2009},
  publisher={SIAM}
}

@article{mathai2002product,
  title={On product moments of order statistics},
  author={Mathai, A.M. and Provost, S.B.},
  journal={Statistical Methods},
  volume={4},
  pages={75--98},
  year={2002}
}

@article{neumann2013central,
  title={A central limit theorem for triangular arrays of weakly dependent random variables, with applications in statistics},
  author={Neumann, M.H.},
  journal={ESAIM: Probability and Statistics},
  volume={17},
  pages={120--134},
  year={2013},
  publisher={EDP Sciences}
}

@article{duffy2011sample,
  title={Sample path large deviations for order statistics},
  author={Duffy, Ken R. and Macci, Claudio and Torrisi, Giovanni Luca},
  journal={Journal of applied probability},
  volume={48},
  number={1},
  pages={238--257},
  year={2011},
  publisher={Cambridge University Press}
}

@article{derrida1998exactly,
  title={An exactly soluble non-equilibrium system: the asymmetric simple exclusion process},
  author={Derrida, Bernard},
  journal={Physics Reports},
  volume={301},
  number={1-3},
  pages={65--83},
  year={1998},
  publisher={Elsevier}
}

@article{frassek2022exact,
  title={Exact solution of an integrable non-equilibrium particle system},
  author={Frassek, Rouven and Giardin{\`a}, Cristian},
  journal={Journal of Mathematical Physics},
  volume={63},
  number={10},
  year={2022},
  publisher={AIP Publishing}
}

@article{frassek2020non,
  title={Non-compact quantum spin chains as integrable stochastic particle processes},
  author={Frassek, Rouven and Giardin{\`a}, Cristian and Kurchan, Jorge},
  journal={Journal of Statistical Physics},
  volume={180},
  number={1},
  pages={135--171},
  year={2020},
  publisher={Springer}
}

@article{carinci2025large,
  title={Large deviations and additivity principle for the open harmonic process},
  author={Carinci, Gioia and Franceschini, Chiara and Frassek, Rouven and Giardin{\`a}, Cristian and Redig, Frank},
  journal={Communications in Mathematical Physics},
  volume={406},
  number={5},
  pages={103},
  year={2025},
  publisher={Springer}
}

@article{giardina2025intertwining,
  title={Intertwining and propagation of mixtures for generalized KMP models and harmonic models},
  author={Giardin{\`a}, Cristian and Redig, Frank and van Tol, Berend},
  journal={Journal of Statistical Physics},
  volume={192},
  number={2},
  pages={21},
  year={2025},
  publisher={Springer}
}

@article{bertini2007stochastic,
  title={Stochastic interacting particle systems out of equilibrium},
  author={Bertini, Lorenzo and De Sole, Alberto and Gabrielli, Davide and Jona-Lasinio, Giovanni and Landim, Claudio},
  journal={Journal of Statistical Mechanics: Theory and Experiment},
  volume={2007},
  number={07},
  pages={P07014},
  year={2007},
  publisher={IOP Publishing}
}

@article{carinci2024solvable,
  title={Solvable stationary non equilibrium states},
  author={Carinci, Gioia and Franceschini, Chiara and Gabrielli, Davide and Giardin{\`a}, Cristian and Tsagkarogiannis, Dimitrios},
  journal={Journal of Statistical Physics},
  volume={191},
  number={1},
  pages={10},
  year={2024},
  publisher={Springer}
}

@article{de2024hidden,
  title={Hidden temperature in the KMP model},
  author={De Masi, Anna and Ferrari, Pablo A and Gabrielli, Davide},
  journal={Journal of Statistical Physics},
  volume={191},
  number={11},
  pages={150},
  year={2024},
  publisher={Springer}
}

@article{franceschini2023integrable,
  title={Integrable heat conduction model},
  author={Franceschini, Chiara and Frassek, Rouven and Giardina, Cristian},
  journal={Journal of Mathematical Physics},
  volume={64},
  number={4},
  year={2023},
  publisher={AIP Publishing}
}
\end{document}